%% file: osqp.tex
\pdfoutput=1
\newif\ifpreprint%
\preprinttrue   

\ifpreprint
\documentclass[12pt]{article}

\usepackage{abstract} 

\usepackage{titlesec}
\def\sectionfont{\sffamily\Large\bfseries\boldmath}
\def\subsectionfont{\sffamily\large\bfseries\boldmath}
\def\paragraphfont{\sffamily\normalsize\bfseries\boldmath}
\titleformat*{\section}{\sectionfont}
\titleformat*{\subsection}{\subsectionfont}
\titleformat*{\subsubsection}{\paragraphfont}
\titleformat*{\paragraph}{\paragraphfont}
\titleformat*{\subparagraph}{\paragraphfont}
\usepackage[small,labelfont={bf,sf}]{caption}  
\usepackage{fullpage,graphicx,psfrag,url}
\usepackage{multirow}
\usepackage{enumitem}
\setlist{nolistsep}

\newcounter{algorithmctr}[section]
\renewcommand{\thealgorithmctr}{\thesection.\arabic{algorithmctr}}

\else

\RequirePackage{fix-cm}

\documentclass[smallextended]{svjour3}       

\smartqed  

\journalname{Mathematical Programming Computation}

\makeatletter
\renewcommand\appendix{\par\small
  \setcounter{section}{0}%
  \setcounter{subsection}{0}%
  \def\paragraph{\@startsection{paragraph}{4}{\z@}%
    {-13pt plus-8pt minus-4pt}{\z@}{\small\itshape}}
  \renewcommand\thesection{\@Alph\c@section}}
\makeatother

\fi

\usepackage[english]{babel}
\usepackage{algorithm}
\usepackage[noend]{algpseudocode}
\usepackage{mathtools}

\usepackage{booktabs}           
\usepackage{csvsimple}	
\usepackage{multirow}
\usepackage{adjustbox} 

\usepackage{xcolor}
\definecolor{borderBoxTodo}{rgb}{0,0.309,0.6}
\definecolor{backgroundBoxTodo}{rgb}{0.701,0.854,1}
\definecolor{lineColorTodo}{rgb}{0.701,0.854,1}
\usepackage[colorinlistoftodos,
            backgroundcolor=backgroundBoxTodo,
            bordercolor = borderBoxTodo,
            linecolor = lineColorTodo]{todonotes}
\usepackage{xcolor}

\usepackage{hyperref}
\hypersetup{%
    colorlinks=true, 
    linktoc=all,     
    linkcolor=black,  
    citecolor= black,
    urlcolor=blue
}

\usepackage[
tight-spacing=true,
round-mode=places,
table-number-alignment=left,
table-figures-integer=1,
round-precision=3,
]{siunitx}
\usepackage{numprint}
\npdecimalsign{.}
\nprounddigits{2}
\usepackage{tikz} 
\usepackage{pgfplots}
\usepgfplotslibrary{groupplots}
\pgfplotsset{compat=newest} 
\usepgfplotslibrary{units} 
\usetikzlibrary{matrix} 
\pgfplotsset{%
    filter discard warning=false,
    /pgfplots/ybar legend/.style={
    /pgfplots/legend image code/.code={%
       \draw[##1,/tikz/.cd,
       yshift=-0.3em
       ]
        (0cm,0cm) rectangle (5pt,0.8em);},
   },
    discard if not/.style 2 args={
    x filter/.code={
    \edef\tempa{\thisrow{#1}}
    \edef\tempb{#2}
    \ifx\tempa\tempb
    \else
	
    \fi
}
}
}

\definecolor{osqp}{RGB}{55,126,184}
\newcommand{\osqpmarker}{o}
\newcommand{\osqplinestyle}{solid}
\definecolor{osqpnows}{RGB}{166,206,227}
\definecolor{gurobi}{RGB}{228,26,28}
\newcommand{\gurobimarker}{square}
\newcommand{\gurobilinestyle}{dashed}
\definecolor{mosek}{RGB}{152,78,163}
\newcommand{\moseklinestyle}{dashdotted}
\definecolor{ecos}{RGB}{166,86,40}

\definecolor{qpoases}{RGB}{77,175,74}
\newcommand{\qpoaseslinestyle}{loosely dashdotdotted}

\definecolor{randomqp}{RGB}{55,126,184}

\definecolor{eqqp}{RGB}{228,26,28}

\definecolor{portfolio}{RGB}{152,78,163}

\definecolor{lasso}{RGB}{166,86,40}

\definecolor{huber}{RGB}{77,175,74}

\definecolor{control}{RGB}{255,127,0}

\definecolor{svm}{RGB}{247,129,191}

\newcommand{\nodec}[1]{\multicolumn{1}{r}{\num[zero-decimal-to-integer=true]{#1}}}

\newcommand{\osqpstatstable}[1]
{\begin{tabular}{rr
S[table-figures-integer=4, round-precision=2]
S[table-figures-integer=4, round-precision=2]
}
\toprule
& & {Median} & {Max} \\
\midrule
\csvreader[head to column names, late after line=\\]{./data/results/#1/ratio_setup_solve.csv}{
	median_ratio=\medianratio,
	max_ratio=\maxratio
	}
{\multirow{2}*{Setup/solve time [$\%$]} & Low accuracy & \medianratio & \maxratio}
\csvreader[head to column names, late after line=\\]{./data/results/#1_high_accuracy/ratio_setup_solve.csv}{
	median_ratio=\medianratio,
	max_ratio=\maxratio
}{ & High accuracy & \medianratio & \maxratio}
\midrule
\csvreader[head to column names, late after line=\\]{./data/results/#1/polish_statistics.csv}{
median_time_increase=\time,
max_time_increase=\maxtime,
}
{\multirow{2}*{Polish time increase [$\%$]} & Low accuracy & \time & \maxtime}
\csvreader[head to column names, late after line=\\]{./data/results/#1_high_accuracy/polish_statistics.csv}{
median_time_increase=\time,
max_time_increase=\maxtime}
{ & High accuracy & \time & \maxtime}
\midrule
\csvreader[head to column names, late after line=\\]{./data/results/#1/rho_updates.csv}{median_rho_updates=\medianrho, max_rho_updates=\maxrho}
{\multirow{2}*{Number of $\rho$ updates} & Low accuracy & \nodec{\medianrho} & \nodec{\maxrho}}
\csvreader[head to column names, late after line=\\]{./data/results/#1_high_accuracy/rho_updates.csv}{
	median_rho_updates=\medianrho,
	max_rho_updates=\maxrho}
	{ & High accuracy & \nodec{\medianrho} & \nodec{\maxrho}}
\bottomrule
\\
	  & & {Mean} & \\
	  \midrule
\csvreader[head to column names, late after line=\\]{./data/results/#1/polish_statistics.csv}{percentage_of_success=\success}
	{\multirow{2}*{Polish success [\%]} & Low accuracy & \success &}
\csvreader[head to column names, late after line=\\]{./data/results/#1_high_accuracy/polish_statistics.csv}{percentage_of_success=\success}
	{ & High accuracy & \success & }
\bottomrule
\end{tabular}
}

\newcommand{\tablecompareosqpvsip}[1]
{\begin{tabular}{@{}rr
    S[table-figures-integer=2]
    S[table-figures-integer=2]
    S[table-figures-integer=2]@{}
}
\toprule
& & ${\rm OSQP}$ & ${\rm GUROBI}$ & ${\rm MOSEK}$\\
\midrule
\csvreader[head to column names, late after line=\\]{./data/results/#1/geom_mean.csv}{
	OSQP=\OSQP,
	GUROBI=\GUROBI,
	MOSEK=\MOSEK}
{\multirow{2}*{\shortstack{Shifted geometric\\means}} & Low accuracy & \OSQP & \GUROBI & \MOSEK}
\csvreader[head to column names, late after line=\\]{./data/results/#1_high_accuracy/geom_mean.csv}{
	OSQP_high=\OSQP,
	GUROBI_high=\GUROBI,
	MOSEK_high=\MOSEK
}{ & High accuracy & \OSQP & \GUROBI & \MOSEK}
\midrule
\csvreader[head to column names, late after line=\\]{./data/results/#1/failure_rates.csv}{
	OSQP=\OSQP,
	GUROBI=\GUROBI,
	MOSEK=\MOSEK}
{\multirow{2}*{Failure rates [$\%$]} & Low accuracy & \OSQP & \GUROBI & \MOSEK}
\csvreader[head to column names, late after line=\\]{./data/results/#1_high_accuracy/failure_rates.csv}{
	OSQP_high=\OSQP,
	GUROBI_high=\GUROBI,
	MOSEK_high=\MOSEK
}{ & High accuracy & \OSQP & \GUROBI & \MOSEK}
\bottomrule
\end{tabular}
}

\newcommand{\tablecompareosqpvsall}[1]
{\begin{tabular}{@{}rr
    S[table-figures-integer=2]
    S[table-figures-integer=2]
    S[table-figures-integer=2]
    S[table-figures-integer=2]
    S[table-figures-integer=3]@{}
}
\toprule
& & ${\rm OSQP}$ & ${\rm GUROBI}$ & ${\rm MOSEK}$ & ${\rm ECOS}$ & ${\rm qpOASES}$\\
\midrule
\csvreader[head to column names, late after line=\\]{./data/results/#1/geom_mean.csv}{
	OSQP=\OSQP,
	GUROBI=\GUROBI,
	MOSEK=\MOSEK,
	ECOS=\ECOS,
	qpOASES=\qpOASES
}
{\multirow{2}*{\shortstack{Shifted geometric\\means}} & Low accuracy & \OSQP & \GUROBI & \MOSEK & \ECOS & \qpOASES}
\csvreader[head to column names, late after line=\\]{./data/results/#1_high_accuracy/geom_mean.csv}{
	OSQP_high=\OSQP,
	GUROBI_high=\GUROBI,
	MOSEK_high=\MOSEK,
	ECOS_high=\ECOS,
	qpOASES=\qpOASES
}{ & High accuracy & \OSQP & \GUROBI & \MOSEK & \ECOS & \qpOASES}
\midrule
\csvreader[head to column names, late after line=\\]{./data/results/#1/failure_rates.csv}{
	OSQP=\OSQP,
	GUROBI=\GUROBI,
	MOSEK=\MOSEK,
	ECOS=\ECOS,
	qpOASES=\qpOASES
}
{\multirow{2}*{Failure rates [$\%$]} & Low accuracy & \OSQP & \GUROBI & \MOSEK & \ECOS & \qpOASES}
\csvreader[head to column names, late after line=\\]{./data/results/#1_high_accuracy/failure_rates.csv}{
	OSQP_high=\OSQP,
	GUROBI_high=\GUROBI,
	MOSEK_high=\MOSEK,
	ECOS_high=\ECOS,
	qpOASES=\qpOASES
}{ & High accuracy & \OSQP & \GUROBI & \MOSEK & \ECOS & \qpOASES}
\bottomrule
\end{tabular}
}

\usepackage[xindy, acronyms]{glossaries}   

\input defs.tex

\ifpreprint
\else
\fi

\ifpreprint
\title{\bfseries\sffamily OSQP: An Operator Splitting Solver for Quadratic Programs}
\author{Bartolomeo Stellato, Goran Banjac, Paul Goulart, \\Alberto Bemporad, and Stephen Boyd}
\else
\title{OSQP: An Operator Splitting Solver for Quadratic Programs\thanks{This work was supported by the People Programme (Marie Curie Actions) of the European Union Seventh Framework Programme (FP7/2007--2013) under REA grant agreement no 607957 (TEMPO).}}
\author{Bartolomeo Stellato \and
Goran Banjac \and
Paul Goulart \and
Alberto Bemporad \and
Stephen Boyd}

\institute{
B. Stellato \at
Massachusetts Institute of Technology, 77 Massachusetts Avenue, Cambridge, MA 02139, USA\\
\email{stellato@mit.edu}           
\and
G. Banjac \at
ETH Z\"{u}rich, R\"{a}mistrasse 101, 8092 Z\"{u}rich, Switzerland\\
\email{gbanjac@control.ee.ethz.ch}           
\and
P. Goulart \at
University of Oxford, Parks Road, Oxford, OX1 3PJ, UK\\
\email{paul.goulart@eng.ox.ac.uk}           
\and
A. Bemporad \at
IMT Institute for Advanced Studies Lucca, Piazza S.\ Francesco 19, 55100 Lucca, Italy\\
\email{alberto.bemporad@imtlucca.it}           
\and
S. Boyd \at
Stanford University, 450 Serra Mall, Stanford CA 94305, USA\\
\email{boyd@stanford.edu}           
}

\date{Received: date / Accepted: date}

\fi

\begin{document}
\maketitle

\begin{abstract}
We present a general-purpose solver for convex quadratic programs based on the \glsdesc{ADMM}, employing a novel operator splitting technique that requires the solution of a quasi-definite linear system with the same coefficient matrix at almost every iteration.
Our algorithm is very robust, placing no requirements on the problem data such as positive definiteness of the objective function or linear independence of the constraint functions.
It can be configured to be division-free once an initial matrix factorization is carried out, making it suitable for real-time applications in embedded systems.
In addition, our technique is the first operator splitting method for quadratic programs able to reliably detect primal and dual infeasible problems from the algorithm iterates.
The method also supports factorization caching and warm starting, making it particularly efficient when solving parametrized problems arising in finance, control, and machine learning.
Our open-source C implementation OSQP has a small footprint, is library-free, and has been extensively tested on many problem instances from a wide variety of application areas.
It is typically ten times faster than competing interior-point methods, and sometimes much more when factorization caching or warm start is used.
OSQP has already shown a large impact with tens of thousands of users both in academia and in large corporations.

\ifpreprint \else
\fi

\end{abstract}

\input{sections/introduction.tex}

\input{sections/optimality_conditions.tex}

\input{sections/certificates_infeasibility.tex}

\input{sections/solution.tex}

\input{sections/polishing.tex}

\input{sections/prec_param_selection.tex}

\input{sections/parametric_programs.tex}

\input{sections/osqp_implementation.tex}

\input{sections/examples.tex}

\input{sections/conclusion.tex}

\input{sections/appendix.tex}

\bibliography{refs}

\end{document}

%% file: defs.tex
\newcommand{\eg}{{\it e.g.}}
\newcommand{\ie}{{\it i.e.}}

\newcommand{\ones}{\mathbf 1}
\newcommand{\reals}{{\mbox{\bf R}}}

\newcommand{\symm}{{\mbox{\bf S}}}  

\newcommand{\diag}{\mathop{\bf diag}}

\newcommand{\tpose}{T}
\newcommand{\eps}{\varepsilon}
\newcommand{\eqdef}{=}

\newcommand{\Cpp}{{C\nolinebreak[4]\hspace{-.05em}\raisebox{.4ex}{\tiny\bf ++}}}

\newcommand{\Chap}{{Chap.\ }}
\newcommand{\Prop}{{Prop.\ }}

\newcommand{\Thm}{{Thm.\ }}
\newcommand{\Lem}{{Lem.\ }}

\newcommand{\mean}{\mathop{\rm mean}}
\newcommand{\argmin}{\mathop{\rm argmin}}

\newcommand{\nnz}{\mathop{\rm nnz{}}}

\newacronym{LP}{LP}{linear program}
\newacronym{QP}{QP}{quadratic program}
\newacronym{MIQP}{MIQP}{mixed-integer quadratic program}
\newacronym{MPC}{MPC}{model predictive control}
\newacronym{MHE}{MHE}{moving horizon estimation}
\newacronym{SQP}{SQP}{sequential quadratic programming}
\newacronym{ADMM}{ADMM}{alternating direction method of multipliers}

%% file: sections/introduction.tex

\section{Introduction}

\subsection{The problem}
Consider the following optimization problem
\begin{equation}\label{pb:main}
\begin{array}{ll}
	\mbox{minimize}   & (1/2) x^\tpose P x + q^\tpose x \\
	\mbox{subject to} & A x \in \mathcal{C},
\end{array}
\end{equation}
where $x\in \reals^{n}$ is the decision variable.
The objective function is defined by a positive semidefinite matrix $P \in \symm_{+}^{n}$ and a vector $q \in \reals^{n}$, and the constraints by a matrix $A\in \reals^{m \times n}$ and a nonempty, closed and convex set $\mathcal{C} \subseteq \reals^m$.
We will refer to it as {\em general (convex) \glsdesc{QP}}.

If the set $\mathcal{C}$ takes the form
\begin{equation*}
	\mathcal{C} = [l,u] \eqdef \left\{z \in \reals^{m} \mid l_i\le z_i \le u_i, \,\, i=1,\dots,m\right\},
\end{equation*}
with ${l_i \in \lbrace-\infty\rbrace\cup\reals}$ and ${u_i \in \reals\cup\lbrace+\infty\rbrace}$, we can write problem~\eqref{pb:main} as
\begin{equation}\label{pb:qp}
\begin{array}{ll}
	\mbox{minimize}   & (1/2) x^\tpose P x + q^\tpose x \\
	\mbox{subject to} & l \leq A x \leq u,
\end{array}
\end{equation}
which we will refer to as a {\em \glsdisp{QP}{quadratic program {\rm (QP)}}}.
Linear equality constraints can be encoded in this way by setting $l_i = u_i$ for some or all of the elements in $(l,u)$.
Note that any \gls{LP} can be written in this form by setting $P=0$.
We will characterize the size of~\eqref{pb:qp} with the tuple $(n, m, N)$ where $N$ is the sum of the number of nonzero entries in $P$ and $A$, \ie, $N \eqdef \nnz(P) + \nnz(A)$.

\paragraph{Applications.} Optimization problems of the form~\eqref{pb:main} arise in a huge variety of applications in engineering, finance, operations research and many other fields.
Applications in machine learning include support vector machines (SVM)~\cite{Cortes:1995}, Lasso~\cite{Tibshirani:1996,Candes:2008} and Huber fitting~\cite{Huber:1964,Huber:1981}.
Financial applications of~\eqref{pb:main} include portfolio optimization~\cite{cornuejols2006,JOFI:JOFI1525,OPT-001,OPT-023}~\cite[\S 4.4.1]{boyd2004convex}.
In the field of control engineering,
\gls{MPC}~\cite{Rawlings:2009,GARCIA1989335} and \gls{MHE}~\cite{Allgower1999} techniques require the solution of a \gls{QP} at each time instant.
Several signal processing problems also fall into the same class~\cite[\S 6.3.3]{boyd2004convex}\cite{5447065}.
In addition, the numerical solution of \gls{QP} subproblems is an essential component in nonconvex optimization methods such as \gls{SQP}~\cite[\Chap 18]{NW06} and mixed-integer optimization using branch-and-bound algorithms~\cite{Belotti:2013dva,FL99}.

\subsection{Solution methods}
Convex \glspl{QP} have been studied since the 1950s~\cite{MargueriteWolfe1956}, following from the seminal work on LPs started by Kantorovich~\cite{Kantorovich:1939}.
Several solution methods for both \glspl{LP} and \glspl{QP} have been proposed and improved upon throughout the years.

\paragraph{Active-set methods.}
Active-set methods were the first algorithms popularized as solution methods for \glspl{QP} \cite{Wolfe1959}, and were obtained from an extension of Dantzig's simplex method for solving \glspl{LP}~\cite{dantzig1963}.
Active-set algorithms select an active-set (\ie, a set of binding constraints) and then iteratively adapt it by adding and dropping constraints from the index of active ones~\cite[\S 16.5]{NW06}.
New active constraints are added based on the cost function gradient and the current dual variables.
Active-set methods for \glspl{QP} differ from the simplex method for \glspl{LP} because the iterates are not necessarily vertices of the feasible region.
These methods can easily be warm started to reduce the number of active-set recalculations required.
However, the major drawback of active-set methods is that the worst-case complexity grows exponentially with the number of constraints, since it may be necessary to investigate all possible active-sets before reaching the optimal one~\cite{KleeMinty:1970}.
Modern implementations of active-set methods for the solution of \glspl{QP} can be found in many commercial solvers, such as MOSEK~\cite{mosek} and GUROBI~\cite{gurobi}, and in the open-source solver qpOASES~\cite{Ferreau:2014}.

\paragraph{Interior-point methods.}
Interior-point algorithms gained popularity in the 1980s as a method for solving \glspl{LP} in polynomial time~\cite{Karmarkar1984,Gill1986}.
In the 90s these techniques were extended to general convex optimization problems, including \glspl{QP}~\cite{doi:10.1137/1.9781611970791}.
Interior-point methods model the problem constraints as parametrized penalty functions, also referred to as \emph{barrier functions}.
At each iteration an unconstrained optimization problem is solved for varying barrier function parameters until the optimum is achieved; see~\cite[\Chap 11]{boyd2004convex} and~\cite[\S 16.6]{NW06} for details.
Primal-dual interior-point methods, in particular the Mehrotra predictor-corrector~\cite{doi:10.1137/0802028} method, became the algorithms of choice for practical implementation~\cite{doi:10.1137/1.9781611971453} because of their good performance across a wide range of problems.
However, interior-point methods are not easily warm started and do not scale well for very large problems.
Interior-point methods are currently the default algorithms in the commercial solvers MOSEK~\cite{mosek}, GUROBI~\cite{gurobi} and CVXGEN~\cite{Mattingley:2012} and in the open-source solver OOQP~\cite{Gertz:2003:OSQ:641876.641880}.

\paragraph{First-order methods.}
First-order optimization methods for solving quadratic programs date to the 1950s~\cite{MargueriteWolfe1956}.
These methods iteratively compute an optimal solution using only first-order information about the cost function.
Operator splitting techniques such as the Douglas-Rachford splitting~\cite{doi:10.1137/0716071,10.2307/1993056} are a particular class of first-order methods which model the optimization problem as the problem of finding a zero of the sum of monotone operators.

In recent years, the operator splitting method known as the \gls{ADMM}~\cite{GABAY197617,Glowinski1975} has received particular attention because of its very good practical convergence behavior; see~\cite{boyd2011distributed} for a survey.
\gls{ADMM} can be seen as a variant of the classical alternating projections algorithm~\cite{doi:10.1137/S0036144593251710} for finding a point in the intersection of two convex sets, and can also be shown to be equivalent to the Douglas-Rachford splitting~\cite{GABAY1983299}.
\gls{ADMM} has been shown to reliably provide modest accuracy solutions to \glspl{QP} in a relatively small number of computationally inexpensive iterations.
It is therefore well suited to applications such as embedded optimization or large-scale optimization, wherein high accuracy solutions are typically not required due to noise in the data and arbitrariness of the cost function.
\gls{ADMM} steps are computationally very cheap and simple to implement, and thus ideal for embedded processors with limited computing resources such as those found in embedded control systems~\cite{6882832,6422363,SYS-008}.
\gls{ADMM} is also compatible with distributed optimization architectures enabling the solution of very large-scale problems~\cite{boyd2011distributed}.

A drawback of first-order methods is that they are typically unable to detect primal and/or dual infeasibility.
In order to address this shortcoming, a homogeneous self-dual embedding has been proposed in conjunction with \gls{ADMM} for solving conic optimization problems and implemented in the open-source solver SCS~\cite{ocpb:16}.
Although every \gls{QP} can be reformulated as a conic program, this reformulation is not efficient from a computational point of view.
A further drawback of \gls{ADMM} is that number of iterations required to converge is highly dependent on the problem data and on the user's choice of the algorithm's step-size parameters.
Despite some recent theoretical results~\cite{7465685,BanjacLinearConvergence}, it remains unclear how to select those parameters to optimize the algorithm convergence rate.
For this reason, even though there are several benefits in using \gls{ADMM} techniques for solving optimization problems, there exists no reliable general-purpose \gls{QP} solver based on operator splitting methods.

\subsection{Our approach}
In this work we present a new general-purpose \gls{QP} solver based on \gls{ADMM} that is able to provide high accuracy solutions.
The proposed algorithm is based on a novel splitting requiring the solution of a quasi-definite linear system that is always solvable for any choice of problem data.
We therefore impose no constraints such as strict convexity of the cost function or linear independence of the constraints.
Since the linear system's matrix coefficients remain the same at every iteration when $\rho$ is fixed, our algorithm requires only a single factorization to solve the QP \eqref{pb:qp}.
Once this initial factorization is computed, we can fix the linear system matrix coefficients to make the algorithm division-free.
If we allow divisions, then we can make occasional updates to the term $\rho$ in this linear system to improve our algorithm's convergence.
We find  that our algorithm typically updates these coefficients very few times, \eg, $1$ or $2$ in our experiments.
In contrast to other first-order methods, our approach is able to return primal and dual solutions when the problem is solvable or to provide certificates of primal and dual infeasibility without resorting to the homogeneous self-dual embedding.

To obtain high accuracy solutions, we perform \emph{solution polishing} on the iterates obtained from \gls{ADMM}.
By identifying the active constraints from the final dual variable iterates, we construct an ancillary equality-constrained \gls{QP} whose solution is equivalent to that of the original \gls{QP}~\eqref{pb:main}.
This ancillary problem is then solved by computing the solution of a single linear system of typically much lower dimensions than the one solved during the \gls{ADMM} iterations.
If we identify the active constraints correctly, then the resulting solution of our method has accuracy equal to or even better than interior-point methods.

Our algorithm can be efficiently warm started to reduce the number of iterations.
Moreover, if the problem matrices do not change then the quasi-definite system factorization can be reused across multiple solves greatly improving the computation time.
This feature is particularly useful when solving multiple instances of parametric \glspl{QP} where only a few elements of the problem data change.
Examples illustrating the effectiveness of the proposed algorithm in parametric programs arising in embedded applications appear in~\cite{embedded_osqp:2017}.

We implemented our method in the open-source ``Operator Splitting Quadratic Program'' (OSQP) solver.
OSQP is written in C and can be compiled to be library free.
OSQP is robust against noisy and unreliable problem data, has a very small code footprint, and is suitable for both embedded and large-scale applications.
We have extensively tested our code and carefully tuned its parameters by solving millions of \glspl{QP}.
We benchmarked our solver against state-of-the-art interior-point and active-set solvers over a benchmark library of 1400 problems from 7 different classes and over the hard \glspl{QP} Maros-M\'{e}sz\'{a}ros test set~\cite{maros1999}.
Numerical results show that our algorithm is able to provide up to an order of magnitude computational time improvements over existing commercial and open-source solvers in a wide variety of applications.
We also showed further time reductions from warm starting and factorization caching.

%% file: sections/optimality_conditions.tex

\section{Optimality conditions}

We will find it convenient to rewrite problem~\eqref{pb:main} by introducing an additional decision variable $z\in\reals^m$, to obtain the equivalent problem
\begin{equation}\label{pb:equivalent}
\begin{array}{ll}
	\mbox{minimize}   & (1/2) x^\tpose P x + q^\tpose x \\
	\mbox{subject to} & A x = z \\
	& z \in \mathcal{C}.
\end{array}
\end{equation}
We can write the optimality conditions of problem~\eqref{pb:equivalent} as~\cite[\Lem A.1]{osqp-convergence}~\cite[\Thm 6.12]{rockafellar1998variational}
\begin{align}
    & Ax = z, \label{eq:prim_feas} \\
    & Px + q + A^\tpose y = 0, \label{eq:dual_feas}\\
    & z \in \mathcal{C}, \quad y\in N_{\mathcal{C}}(z),\label{eq:separating_hyperplane}
\end{align}
where $y\in \reals^{m}$ is the Lagrange multiplier associated with the constraint $Ax=z$ and $N_{\mathcal{C}}(z)$ denotes the normal cone of $\mathcal{C}$ at $z$.
If there exist $x \in \reals^n$, $z \in \reals^{m}$ and $y \in \reals^{m}$ that satisfy the conditions above, then we say that $(x,z)$ is a {\em primal} and $y$ is a {\em dual} solution to problem~\eqref{pb:equivalent}.
We define the primal and dual residuals of problem~\eqref{pb:main} as
\begin{align}
    r_{\rm prim} &\eqdef Ax - z, \label{eq:res_prim}\\
    r_{\rm dual} &\eqdef Px + q + A^\tpose y. \label{eq:res_dual}
\end{align}

\paragraph{Quadratic programs.}
In case of \glspl{QP} of the form~\eqref{pb:qp}, condition \eqref{eq:separating_hyperplane} reduces to
\begin{equation}\label{eq:compl_slack}
    l \le z \le u, \qquad y_{+}^\tpose (z - u) = 0,\qquad y_{-}^\tpose (z - l) = 0,
\end{equation}
where $y_{+} \eqdef \max(y, 0)$ and $y_{-} \eqdef \min(y, 0)$.



%% file: sections/certificates_infeasibility.tex

\subsection{Certificates of primal and dual infeasibility}

From the \emph{theorem of strong alternatives}~\cite[\S 5.8]{boyd2004convex},~\cite[\Prop 3.1]{osqp-convergence}, exactly one of the following sets is nonempty
\begin{align}
	\mathcal{P} &= \left\{x \in \reals^{n}\mid Ax\in\mathcal{C} \right\}, \\
	\mathcal{D} &= \left\{y \in \reals^{m} \mid A^\tpose y = 0, \quad S_\mathcal{C}(y)<0 \right\},\label{eq:infeas-general-D}
\end{align}%
\label{eq:infeas-general}%
where $S_\mathcal{C}$ is the support function of $\mathcal{C}$, provided that some type of constraint qualification holds~\cite{boyd2004convex}.
In other words, any variable $y \in \mathcal{D}$ serves as a \emph{certificate} that problem~\eqref{pb:main} is primal infeasible.

\paragraph{Quadratic programs.}
In case $\mathcal{C}=[l,u]$, certifying primal infeasibility of~\eqref{pb:qp} amounts to finding a vector $y\in\reals^m$ such that
\begin{equation}\label{eq:prim-infeas-linear}
		A^\tpose y = 0, \quad u^\tpose y_{+} + l^\tpose y_{-} < 0.
\end{equation}
Similarly, it can be shown that a vector $x\in\reals^n$ satisfying
\begin{equation}\label{eq:dual-infeas-linear}
		Px = 0, \quad q^\tpose x < 0, \quad (Ax)_i
		\left\{
		\begin{aligned} =0 &\quad l_i,u_i\in\reals \\
		\geq 0 &\quad  u_i = +\infty, \,l_i \in \reals \\
		\leq 0 &\quad l_i = -\infty, \,u_i \in \reals
	\end{aligned}
	\right.
\end{equation}
is a certificate of dual infeasibility for problem~\eqref{pb:qp}; see \cite[\Prop 3.1]{osqp-convergence} for more details.

%% file: sections/solution.tex

\section{Solution with ADMM}
\label{sec:solution_with_admm}
Our method solves problem~\eqref{pb:equivalent} using \gls{ADMM}~\cite{boyd2011distributed}.
By introducing auxiliary variables $\tilde{x}\in\reals^n$ and $\tilde{z}\in\reals^m$, we can rewrite problem \eqref{pb:equivalent} as
\begin{equation}\label{pb:admm_qp}
\begin{array}{ll}
	\mbox{minimize}   & (1/2) \tilde{x}^\tpose P \tilde{x} + q^\tpose \tilde{x} + \mathcal{I}_{Ax=z}(\tilde{x},\tilde{z}) + \mathcal{I}_{\mathcal{C}}(z) \\
	\mbox{subject to} & (\tilde{x}, \tilde{z}) = (x, z),
\end{array}
\end{equation}
where $\mathcal{I}_{Ax=z}$ and $\mathcal{I}_{\mathcal{C}}$ are the indicator functions given by
\[
	\mathcal{I}_{Ax=z}(x,z) = \begin{cases} 0 & Ax=z \\ +\infty & \text{otherwise} \end{cases},
	\qquad
	\mathcal{I}_\mathcal{C}(z) = \begin{cases} 0 & z\in\mathcal{C} \\ +\infty & \text{otherwise} \end{cases}.
\]

An iteration of \gls{ADMM} for solving problem \eqref{pb:admm_qp} consists of the following steps:
\ifpreprint
\begin{align}
	(\tilde{x}^{k+1}, \tilde{z}^{k+1}) &\gets \begin{multlined}[t][.7\textwidth]\argmin_{(\tilde{x},\tilde{z}): A\tilde{x}=\tilde{z}} (1/2) \tilde{x}^\tpose P \tilde{x} + q^\tpose \tilde{x}
	+ (\sigma/2) \| \tilde{x} - x^k + \sigma^{-1} w^{k} \|_2^2  \\
	+ (\rho/2) \| \tilde{z} - z^k + \rho^{-1} y^{k} \|_2^2 \end{multlined}\label{eq:ADMM_iter1} \\
    x^{k+1} &\gets \alpha \tilde{x}^{k+1} + (1-\alpha)x^{k} + \sigma^{-1} w^{k} \label{eq:ADMM_iter2} \\
    z^{k+1} &\gets \Pi\left(\alpha \tilde{z}^{k+1} + (1-\alpha)z^{k} + \rho^{-1} y^{k} \right) \label{eq:ADMM_iter3} \\
    w^{k+1} &\gets w^{k} + \sigma \left( \alpha \tilde{x}^{k+1} + (1-\alpha)x^{k} - x^{k+1} \right) \label{eq:ADMM_iter4} \\
	y^{k+1} &\gets y^{k} + \rho \left( \alpha \tilde{z}^{k+1} + (1-\alpha)z^{k} - z^{k+1} \right) \label{eq:ADMM_iter5}
\end{align}
\else
\begin{align}
 \hspace{-.75em}(\tilde{x}^{k+1}, \tilde{z}^{k+1}) &\gets\hspace{-.5em}\begin{multlined}[t][.7\textwidth] \argmin_{(\tilde{x},\tilde{z}): A\tilde{x}=\tilde{z}} (1/2) \tilde{x}^\tpose P \tilde{x} + q^\tpose \tilde{x} \\
	+ (\sigma/2) \| \tilde{x} - x^k + \sigma^{-1} w^{k} \|_2^2
	+ (\rho/2) \| \tilde{z} - z^k + \rho^{-1} y^{k} \|_2^2\end{multlined} \label{eq:ADMM_iter1} \\
    x^{k+1} &\gets \alpha \tilde{x}^{k+1} + (1-\alpha)x^{k} + \sigma^{-1} w^{k} \label{eq:ADMM_iter2} \\
    z^{k+1} &\gets \Pi\left(\alpha \tilde{z}^{k+1} + (1-\alpha)z^{k} + \rho^{-1} y^{k} \right) \label{eq:ADMM_iter3} \\
    w^{k+1} &\gets w^{k} + \sigma \left( \alpha \tilde{x}^{k+1} + (1-\alpha)x^{k} - x^{k+1} \right) \label{eq:ADMM_iter4} \\
	y^{k+1} &\gets y^{k} + \rho \left( \alpha \tilde{z}^{k+1} + (1-\alpha)z^{k} - z^{k+1} \right) \label{eq:ADMM_iter5}
\end{align}
\fi
where $\sigma>0$ and $\rho>0$ are the \emph{step-size parameters}, $\alpha \in (0,2)$ is the \emph{relaxation parameter}, and $\Pi$ denotes the Euclidean projection onto $\mathcal{C}$.
The introduction of the splitting variable $\tilde{x}$ ensures that the subproblem in~\eqref{eq:ADMM_iter1} is always solvable for any $P\in\symm_+^n$ which can also be $0$ for \glspl{LP}.
Note that all the derivations hold also for $\sigma$ and $\rho$ being positive definite diagonal matrices.
The iterates $w^k$ and $y^k$ are associated with the dual variables of the equality constraints $\tilde{x} = x$ and $\tilde{z} = z$, respectively.
Observe from steps \eqref{eq:ADMM_iter2} and \eqref{eq:ADMM_iter4} that $w^{k+1}=0$ for all $k\ge 0$, and consequently the $w$-iterate and the step \eqref{eq:ADMM_iter4} can be disregarded.

\subsection{Solving the linear system}\label{sec:solving_linear_system}
Evaluating the \gls{ADMM} step~\eqref{eq:ADMM_iter1} involves solving the equality constrained \gls{QP}
\ifpreprint
\begin{equation}\label{pb:admm_step1}
\begin{array}{ll}
\mbox{minimize}   & (1/2) \tilde{x}^\tpose P \tilde{x} +\! q^\tpose \tilde{x} +\! (\sigma/2) \| \tilde{x} - x^k \|_2^2 +\! (\rho/2) \| \tilde{z} - z^k +\! \rho^{-1} y^{k} \|_2^2 \\
\mbox{subject to} & A\tilde{x} = \tilde{z}.
\end{array}
\end{equation}
\else
\begin{equation}\label{pb:admm_step1}
\setlength{\arraycolsep}{2pt}
\begin{array}{ll}
\hspace{-.8em}\mbox{minimize}   & (1/2) \tilde{x}^\tpose P \tilde{x} +\! q^\tpose \tilde{x} +\! (\sigma/2) \| \tilde{x} - x^k \|_2^2 +\! (\rho/2) \| \tilde{z} - z^k +\! \rho^{-1} y^{k} \|_2^2 \\
\hspace{-.8em}\mbox{subject to} & A\tilde{x} = \tilde{z}.
\end{array}
\end{equation}
\fi
The optimality conditions for this equality constrained \gls{QP} are
\begin{align}
	P\tilde{x}^{k+1} + q + \sigma(\tilde{x}^{k+1} - x^k) + A^\tpose \nu^{k+1} &= 0, \label{admm_step1:eq1} \\
	\rho(\tilde{z}^{k+1}-z^k) + y^k - \nu^{k+1} &= 0, \label{admm_step1:eq2} \\
	A\tilde{x}^{k+1} - \tilde{z}^{k+1} &= 0, \label{admm_step1:eq3}
\end{align}
where $\nu^{k+1} \in \reals^m$ is the Lagrange multiplier associated with the constraint $Ax=z$.
By eliminating the variable $\tilde{z}^{k+1}$ from~\eqref{admm_step1:eq2}, the above linear system reduces to
\begin{align}
    &\begin{bmatrix} P + \sigma I & A^\tpose \\ A & -\rho^{-1}I \end{bmatrix} \begin{bmatrix} \tilde{x}^{k+1} \\ \nu^{k+1} \end{bmatrix}= \begin{bmatrix}\sigma x^{k} - q \\ z^{k} - \rho^{-1} y^k \end{bmatrix}, \label{eq:kkt}
\end{align}
with $\tilde{z}^{k+1}$ recoverable as
\begin{equation*}
\tilde{z}^{k+1} = z^k + \rho^{-1}(\nu^{k+1} - y^{k}).
\end{equation*}
We will refer to the coefficient matrix in~\eqref{eq:kkt} as the \emph{KKT matrix}.
This matrix always has full rank thanks to the positive parameters $\sigma$ and $\rho$ introduced in our splitting, so~\eqref{eq:kkt} always has a unique solution for any matrices $P\in\symm_+^n$ and $A\in\reals^{m\times n}$.
In other words, we do not impose any additional assumptions on the problem data such as strong convexity of the objective function or linear independence of the constraints as was done in~\cite{6892987,raghunathan2014cdc,raghunathan2014}. {
}

\paragraph{Direct method.}
A direct method for solving the linear system \eqref{eq:kkt} computes its solution by first factoring the KKT matrix and then performing forward and backward substitution.
Since the KKT matrix remains the same for every iteration of \gls{ADMM}, we only need to perform the factorization once prior to the first iteration and cache the factors so that we can reuse them in subsequent iterations.
This approach is very efficient when the factorization cost is considerably higher than the cost of forward and backward substitutions, so that each iteration is computed quickly.
Note that if $\rho$ or $\sigma$ change, the KKT matrix needs to be factored again.

Our particular choice of splitting results in a KKT matrix that is quasi-definite, \ie, it can be written as a 2-by-2 block-symmetric matrix where the $(1,1)$-block is positive definite, and the $(2,2)$-block is negative definite.
It therefore always has a well defined $LDL^\tpose$ factorization, with~$L$ being a lower triangular matrix with unit diagonal elements and~$D$ a diagonal matrix with nonzero diagonal elements~\cite{vanderbei1995symmetric}.
Note that once the factorization is carried out, computing the solution of~\eqref{eq:kkt} can be made division-free by storing $D^{-1}$ instead of $D$.

When the KKT matrix is sparse and quasi-definite, efficient algorithms can be used for computing a suitable permutation matrix~$P$ for which the factorization of $PKP^\tpose$ results in a sparse factor~$L$~\cite{Amestoy:2004,davis2006direct} without regard for the actual nonzero values appearing in the KKT matrix.
The $LDL^\tpose$ factorization consists of two steps.
In the first step we compute the sparsity pattern of the factor $L$.
This step is referred to as the {\em symbolic factorization} and requires only the sparsity pattern of the KKT matrix.
In the second step, referred to as the {\em numerical factorization}, we determine the values of nonzero elements in~$L$ and~$D$.
Note that we do not need to update the symbolic factorization if the nonzero entries of the KKT matrix change but the sparsity pattern remains the same.

\paragraph{Indirect method.}
With large-scale \glspl{QP}, factoring linear system~\eqref{eq:kkt} might be prohibitive. In these cases it might be more convenient to use an indirect method by solving instead the linear system
\begin{equation*}
	\left(P + \sigma I + \rho A^\tpose A \right)\tilde{x}^{k+1} = \sigma x^{k} - q + A^\tpose (\rho z^{k} - y^{k})
\end{equation*}
obtained by eliminating~$\nu^{k+1}$ from~\eqref{eq:kkt}. We then compute $\tilde{z}^{k+1}$ as $\tilde{z}^{k+1} = A\tilde{x}^{k+1}$.
Note that the coefficient matrix in the above linear system is always positive definite.
The linear system can therefore be solved with an iterative scheme such as the conjugate gradient method~\cite{Golub:1996:MC:248979,NW06}.
When the linear system is solved up to some predefined accuracy, we terminate the method.
We can also warm start the method using the linear system solution at the previous iteration of \gls{ADMM} to speed up its convergence.
In contrast to direct methods, the complexity of indirect methods does not change if we update $\rho$ and $\sigma$ since there is no factorization required. This allows for more updates to take place without any overhead.

\subsection{Final algorithm}\label{sec:algorithm}

By simplifying the \gls{ADMM} iterations according to the previous discussion, we obtain Algorithm~\ref{alg:osqp_algorithm}.
Steps~\ref{alg:update_z_tikde},~\ref{alg:update_x},~\ref{alg:update_z} and~\ref{alg:update_y} of Algorithm~\ref{alg:osqp_algorithm} are very easy to evaluate since they involve only vector addition and subtraction, scalar-vector multiplication and projection onto a box.
Moreover, they are component-wise separable and can be easily parallelized.
The most computationally expensive part is solving the linear system in Step~\ref{alg:solve_lin_sys}, which can be performed as discussed in Section~\ref{sec:solving_linear_system}.

\begin{algorithm}
\caption{}
\label{alg:osqp_algorithm}
\begin{algorithmic}[1]
\State {\bf given} initial values $x^0$, $z^0$, $y^0$ and parameters $\rho>0$, $\sigma>0$, $\alpha \in (0,2)$
\Repeat
\State $(\tilde{x}^{k+1}, \nu^{k+1}) \gets$ solve linear system $\begin{bmatrix} P + \sigma I & A^\tpose \\ A & -\rho^{-1}I \end{bmatrix} \begin{bmatrix} \tilde{x}^{k+1} \\ \nu^{k+1} \end{bmatrix}= \begin{bmatrix}\sigma x^{k} - q \\ z^{k} - \rho^{-1} y^k \end{bmatrix}$ \label{alg:solve_lin_sys}
\State $\tilde{z}^{k+1} \gets z^k + \rho^{-1}(\nu^{k+1} - y^{k})$ \label{alg:update_z_tikde}
\State $x^{k+1} \gets \alpha \tilde{x}^{k+1} + (1-\alpha)x^{k}$ \label{alg:update_x}
\State $z^{k+1} \gets \Pi\left(\alpha \tilde{z}^{k+1} + (1-\alpha)z^{k} + \rho^{-1} y^{k} \right)$ \label{alg:update_z}
\State $y^{k+1} \gets y^{k} + \rho \left( \alpha \tilde{z}^{k+1} + (1-\alpha)z^{k} - z^{k+1} \right)$ \label{alg:update_y}
\Until{termination criterion is satisfied 
}
\end{algorithmic}
\end{algorithm}

\subsection{Convergence and infeasibility detection}\label{sec:Convergence}
We show in this section that the proposed algorithm generates a sequence of iterates $(x^{k},z^{k},y^{k})$ that in the limit satisfy the optimality conditions~\eqref{eq:prim_feas}--\eqref{eq:separating_hyperplane} when problem~\eqref{pb:main} is solvable, or provides a certificate of primal or dual infeasibility otherwise.


If we denote the argument of the projection operator in step~\ref{alg:update_z} of Algorithm~\ref{alg:osqp_algorithm} by $v^{k+1}$, then we can express $z^k$ and $y^k$ as
\begin{equation}\label{zy-update}
	z^{k} = \Pi(v^{k})
	\quad\text{and}\quad
	y^{k} = \rho\left( v^{k} - \Pi(v^{k}) \right).
\end{equation}
Observe from \eqref{zy-update} that iterates $z^{k}$ and $y^{k}$ satisfy optimality condition~\eqref{eq:separating_hyperplane} for all $k>0$ by construction~\cite[\Prop 6.46]{bauschke2011convex}.
Therefore, it only remains to show that optimality conditions \eqref{eq:prim_feas}--\eqref{eq:dual_feas} are satisfied in the limit.

As shown in~\cite[\Prop 5.3]{osqp-convergence}, if problem~\eqref{pb:qp} is solvable, then Algorithm~\ref{alg:osqp_algorithm} produces a convergent sequence of iterates $(x^k,z^k,y^k)$ so that
\begin{align*}
	\lim_{k\to\infty} r_{\rm prim}^k = 0, \\
	\lim_{k\to\infty} r_{\rm dual}^k = 0,
\end{align*}
where $r_{\rm prim}^k$ and $r_{\rm dual}^k$ correspond to the residuals defined in~\eqref{eq:res_prim} and~\eqref{eq:res_dual} respectively.

On the other hand, if problem~\eqref{pb:qp} is primal and/or dual infeasible, then the sequence of iterates $(x^k,z^k,y^k)$ generated by Algorithm~\ref{alg:osqp_algorithm} does not converge.
However, the sequence
\[
	(\delta x^{k}, \delta z^{k}, \delta y^{k} ) \eqdef (x^{k} - x^{k-1}, z^{k} - z^{k-1}, y^{k} - y^{k-1})
\]
always converges and can be used to certify infeasibility of the problem.
According to \cite[\Thm 5.1]{osqp-convergence}, if the problem is primal infeasible, then  $\delta y \eqdef \lim_{k\to\infty} \delta y^{k}$ satisfies conditions~\eqref{eq:prim-infeas-linear}, whereas $\delta x \eqdef \lim_{k\to\infty} \delta x^{k}$ satisfies conditions~\eqref{eq:dual-infeas-linear} if it is dual infeasible.

\subsection{Termination criteria}
\label{sec:termination criteria}

We can define termination criteria for Algorithm~\ref{alg:osqp_algorithm} so that the iterations stop when either a primal-dual solution or a certificate of primal or dual infeasibility is found up to some predefined accuracy.

A reasonable termination criterion for detecting optimality is that the norms of the residuals $r_{\rm prim}^k$ and $r_{\rm dual}^k$ are smaller than some tolerance levels $\eps_{\rm prim}>0$ and $\eps_{\rm dual}>0$ \cite{boyd2011distributed}, \ie,
\begin{equation}\label{eq:convergence_check}
    \| r_{\rm prim}^{k} \|_{\infty} \le \eps_{\rm prim},
    \quad
    \| r_{\rm dual}^{k} \|_{\infty} \le \eps_{\rm dual}.
\end{equation}
We set the tolerance levels as
\begin{align*}
    \eps_{\rm prim} &\eqdef \eps_{\rm abs} + \eps_{\rm rel} \max\lbrace \|Ax^{k}\|_{\infty}, \| z^{k} \|_{\infty} \rbrace \\
    \eps_{\rm dual} &\eqdef \eps_{\rm abs} + \eps_{\rm rel} \max\lbrace \| P x^{k} \|_{\infty}, \| A^\tpose y^{k} \|_{\infty}, \| q \|_{\infty} \rbrace,
\end{align*}
where $\eps_{\rm abs}>0$ and $\eps_{\rm rel}>0$ are absolute and relative tolerances, respectively.

\paragraph{Quadratic programs infeasibility.}
If $\mathcal{C}=[l,u]$, we check the following conditions for primal infeasibility
\begin{equation*}
	\left\|A^\tpose \delta y^{k} \right\|_{\infty} \le \eps_{\rm pinf} \|\delta y^{k}\|_{\infty},
		\quad
	u^\tpose (\delta y^{k})_{+} + l^\tpose (\delta y^{k})_{-} \le \eps_{\rm pinf} \|\delta y^{k}\|_{\infty},
\end{equation*}
where $\eps_{\rm pinf} > 0$ is some tolerance level. Similarly, we define the following criterion for detecting dual infeasibility
\begin{gather*}
		\| P \delta x^{k} \|_{\infty} \le \eps_{\rm dinf} \|\delta x^{k}\|_{\infty}, \quad
		q^\tpose \delta x^{k} \le \eps_{\rm dinf}\|\delta x^{k}\|_{\infty}, \\
		(A \delta x^{k} )_i \begin{cases}
		\in \left[-\eps_{\rm dinf}, \eps_{\rm dinf}\right] \|\delta x^{k}\|_{\infty}& u_i, l_i \in \reals\\
		\ge -\eps_{\rm dinf} \|\delta x^{k}\|_{\infty} & u_i = +\infty\\
	\le \eps_{\rm dinf} \|\delta x^{k}\|_{\infty} &l_i = -\infty,
\end{cases}
\end{gather*}
for $i = 1,\dots,m$ where  $\eps_{\rm dinf} > 0$ is some tolerance level.
Note that $\|\delta x^{k}\|_{\infty}$ and $\|\delta y^{k}\|_{\infty}$ appear in the right-hand sides to avoid division when considering normalized vectors $\delta x^{k}$ and $\delta y^{k}$ in the termination criteria.

%% file: sections/polishing.tex

\section{Solution polishing}
\label{sec:solution_polishing}
Operator splitting methods are typically used for obtaining solution of an optimization problem with a low or medium accuracy.
However, even if a solution is not very accurate we can often guess which constraints are active from an approximate primal-dual solution.
When dealing with \glspl{QP} of the form~\eqref{pb:qp}, we can obtain high accuracy solutions from the final \gls{ADMM} iterates by solving one additional system of equations.

Given a dual solution $y$ of the problem, we define the sets of lower- and upper-active constraints
\begin{align*}
	\mathcal{L} &\eqdef \left\lbrace i \in \{ 1,\dots,m\} \mid y_i < 0 \right\rbrace, \\
	\mathcal{U} &\eqdef \left\lbrace i \in \{ 1,\dots,m\} \mid y_i > 0 \right\rbrace.
\end{align*}
According to~\eqref{eq:compl_slack} we have that $z_\mathcal{L}=l_\mathcal{L}$ and $z_\mathcal{U}=u_\mathcal{U}$, where $l_{\mathcal{L}}$ denotes the vector composed of elements of $l$ corresponding to the indices in $\mathcal{L}$.
Similarly, we will denote by $A_{\mathcal{L}}$ the matrix composed of rows of $A$ corresponding to the indices in $\mathcal{L}$.

If the sets of active constraints are known \emph{a priori}, then a primal-dual solution $(x,y,z)$ can be found by solving the following linear system
\begin{align}\label{eq:kkt_polish}
	&\begin{bmatrix} P & A_{\mathcal{L}}^\tpose & A_{\mathcal{U}}^\tpose \\ A_{\mathcal{L}} & & \\ A_{\mathcal{U}} & & \end{bmatrix}
	\begin{bmatrix} x \\ y_{\mathcal{L}} \\ y_{\mathcal{U}} \end{bmatrix} =
	\begin{bmatrix} -q \\ l_{\mathcal{L}} \\ u_{\mathcal{U}}\end{bmatrix}, \\
	&y_i = 0, \quad i \notin(\mathcal{L}\cup\mathcal{U}), \\
	&z = Ax.
\end{align}

We can then apply the aforementioned procedure to obtain a candidate solution $(x,y,z)$.
If $(x,y,z)$ satisfies the optimality conditions~\eqref{eq:prim_feas}--\eqref{eq:separating_hyperplane}, then our guess is correct and $(x,y,z)$ is a primal-dual solution of problem~\eqref{pb:equivalent}.
This approach is referred to as \emph{solution polishing}.
Note that the dimension of the linear system~\eqref{eq:kkt_polish} is usually much smaller than the KKT system in Section~\ref{sec:solving_linear_system} because the number of active constraints at optimality is less than or equal to $n$ for non-degenerate \glspl{QP}.

However, the linear system~\eqref{eq:kkt_polish} is not necessarily solvable even if the sets of active constraints $\mathcal{L}$ and $\mathcal{U}$ have been correctly identified.
This can happen, \eg, if the solution is degenerate, \ie, if it has one or more redundant active constraints.
We make the solution polishing procedure more robust by solving instead the following linear system
\begin{equation}\label{eq:kkt_polish_reg}
	\begin{bmatrix} P + \delta I & A_{\mathcal{L}}^\tpose & A_{\mathcal{U}}^\tpose \\ A_{\mathcal{L}} & - \delta I & \\\ A_{\mathcal{U}} & & - \delta I \end{bmatrix}
	\begin{bmatrix} \hat{x} \\ \hat{y}_{\mathcal{L}} \\ \hat{y}_{\mathcal{U}} \end{bmatrix} =
	\begin{bmatrix} -q\\ l_{\mathcal{L}} \\ u_{\mathcal{U}} \end{bmatrix},
\end{equation}
where $\delta > 0$ is a regularization parameter with value $\delta \approx 10^{-6}$.
Since the regularized matrix in~\eqref{eq:kkt_polish_reg} is quasi-definite, the linear system~\eqref{eq:kkt_polish_reg} is always solvable.

By using regularization, we actually solve a perturbed linear system and thus introduce a small error to the polished solution.
If we denote by $K$ and $(K+\Delta K)$ the coefficient matrices in~\eqref{eq:kkt_polish} and~\eqref{eq:kkt_polish_reg}, respectively, then we can represent the two linear systems as $Kt = g$ and $(K+\Delta K)\hat{t} = g$.
To compensate for this error, we apply an {\em iterative refinement} procedure~\cite{Wilkinson:1963}, \ie, we iteratively solve
\begin{equation}\label{eq:iterative_refinement}
	(K+\Delta K) \Delta \hat{t}^k = g - K \hat{t}^k
\end{equation}
and update $\hat{t}^{k+1} \eqdef \hat{t}^k + \Delta \hat{t}^k$.
The sequence $\lbrace \hat{t}^k\rbrace$ converges to the true solution $t$, provided that it exists.
Observe that, compared to solving the linear system~\eqref{eq:kkt_polish_reg}, iterative refinement requires only a backward- and a forward-solve, and does not require another matrix factorization.
Since the iterative refinement iterations converge very quickly in practice, we just run them for a fixed number of passes without imposing any termination condition to satisfy. Note that this is the same strategy used in commercial linear system solvers using iterative refinement~\cite{intel2017}.

%% file: sections/prec_param_selection.tex

\section{Preconditioning and parameter selection}
A known weakness of first-order methods is their inability to deal effectively with ill-conditioned problems, and their convergence rate can vary significantly when data are badly scaled.
In this section we describe how to precondition the data and choose the optimal parameters to speed up the convergence of our algorithm.

\subsection{Preconditioning}

Preconditioning is a common heuristic aiming to reduce the number of iterations in first-order methods~\cite[\Chap 5]{NW06},\cite{6892987,BENZI2002418,6126441,GB:15,7465685}.
The optimal choice of preconditioners has been studied for at least two decades and remains an active area of research~\cite[\Chap 2]{doi:10.1137/1.9781611970944},\cite[\Chap 10]{doi:10.1137/1.9781611970937}.
For example, the optimal diagonal preconditioner required to minimize the condition number of a matrix can be found exactly by solving a semidefinite program~\cite{doi:10.1137/1.9781611970777}.
However, this computation is typically \emph{more} complicated than solving the original \gls{QP}, and is therefore unlikely to be worth the effort since preconditioning is only a heuristic to minimize the number of iterations.

In order to keep the preconditioning procedure simple, we instead make use of a simple heuristic called \emph{matrix equilibration}~\cite{bradley2010algorithms,TJ:14,fougnerboyd2017,diamond2017}.
Our goal is to rescale the problem data to reduce the condition number of the symmetric matrix $M\in \symm^{n+m}$ representing the problem data, defined as
\begin{equation}\label{eq:equil_matrix}
	M \eqdef \begin{bmatrix} P & A^\tpose\\ A & 0
	\end{bmatrix}.
\end{equation}
In particular, we use \emph{symmetric matrix equilibration} by computing the diagonal matrix $S \in \symm^{n+m}_{++}$ to decrease the condition number of $S M S$.
We can write matrix~$S$ as
\begin{equation}\label{eq:diagonal_matrix_equil}
	S = \begin{bmatrix} D & \\ & E \end{bmatrix},
\end{equation}
where $D \in \symm^n_{++}$ and $E \in \symm^m_{++}$ are both diagonal.
In addition, we would like to normalize the cost function to prevent the dual variables from being too large.
We can achieve this by multiplying the cost function by the scalar $c > 0$.

Preconditioning effectively modifies problem~\eqref{pb:main} into the following \begin{equation}\label{pb:scaled_qp}
\begin{array}{ll}
	\mbox{minimize}   & (1/2) \bar{x}^\tpose \bar{P} \bar{x} + \bar{q}^\tpose \bar{x}\\[.25em]
	\mbox{subject to} & \bar{A} \bar{x} \in \bar{\mathcal{C}},
\end{array}
\end{equation}
where $\bar{x}=D^{-1}x$, $\bar{P} = cDPD$, $\bar{q} = cDq$, $\bar{A} = E A D$ and $\bar{\mathcal{C}} \eqdef \lbrace Ez \in \reals^m \mid z \in \mathcal{C} \rbrace$.
The dual variables of the new problem are $\bar{y}=cE^{-1}y$.
Note that when $\mathcal{C}=[l,u]$ the Euclidean projection onto $\bar{\mathcal{C}}=[El,Eu]$ is as easy to evaluate as the projection onto~$\mathcal{C}$.

The main idea of the equilibration procedure is to scale the rows of matrix $M$ so that they all have equal $\ell_p$ norm.
It is possible to show that finding such a scaling matrix $S$ can be cast as a convex optimization problem~\cite{1429569}.
However, it is computationally more convenient to solve this problem with heuristic iterative methods, rather than continuous optimization algorithms such as interior-point methods.
We refer the reader to~\cite{bradley2010algorithms} for more details on matrix equilibration.

\paragraph{Ruiz equilibration.}
In this work we apply a variation of the Ruiz equilibration~\cite{ruiz2001}.
This technique was originally proposed to equilibrate square matrices showing fast linear convergence superior to other methods such as the Sinkhorn-Knopp equilibration~\cite{sinkhorn1967concerning}.
Ruiz equilibration converges in few tens of iterations even in cases when Sinkhorn-Knopp equilibration takes thousands of iterations~\cite{knightruizucar2014}.
The steps are outlined in Algorithm~\ref{alg:ruiz_equil} and differ from the original Ruiz algorithm by adding a cost scaling step that takes into account very large values of the cost.
\begin{algorithm}[t]
\caption{Modified Ruiz equilibration}\label{alg:ruiz_equil}
\begin{algorithmic}
\State{\bf initialize $c=1$, $S=I$, $\delta = 0$, $\bar{P}=P, \bar{q}=q, \bar{A}=A, \bar{\mathcal{C}} = \mathcal{C}$}
\While{$\|1 - \delta\|_{\infty} > \eps_{\rm equil}$}
\For{$i = 1,\dots,n+m$}
		\State $\delta_i \gets 1 / \sqrt{\|M_{i}\|_{\infty}}$ \Comment{$M$ equilibration}
\EndFor
	\State $\bar{P}, \bar{q}, \bar{A}, \bar{\mathcal{C}} \gets$ Scale $\bar{P}, \bar{q}, \bar{A}, \bar{\mathcal{C}}$ using $\diag(\delta)$
		\State $\gamma \gets 1/\max\lbrace\mean(\|\bar{P}_{i}\|_{\infty}), \|\bar{q}\|_{\infty}\rbrace$ \Comment{Cost scaling}
		\State $\bar{P} \gets \gamma \bar{P}, \bar{q} \gets \gamma \bar{q}$
	\State $S \gets \diag(\delta) S, c \gets \gamma c$
\EndWhile
\Return{$S, c$}
\end{algorithmic}
\end{algorithm}
The first part is the usual Ruiz equilibration step.
Since $M$ is symmetric, we focus only on the columns $M_{i}$ and apply the scaling to both sides of $M$.
At each iteration, we compute the $\infty$-norm of each column and normalize that column by the inverse of its square root.
The second part is a cost scaling step.
The scalar $\gamma$ is the current cost normalization coefficient taking into account the maximum between the average norm of the columns of $\bar{P}$ and the norm of $\bar{q}$.
We normalize problem data $\bar{P}$, $\bar{q}$, $\bar{A}$, $\bar{l}$, $\bar{u}$ in place at each iteration using the current values of $\delta$ and $\gamma$.

\paragraph{Unscaled termination criteria.}
Although we rescale our problem in the form~\eqref{pb:scaled_qp}, we would still like to apply the stopping criteria defined in Section~\ref{sec:termination criteria} to an unscaled version of our problem.
The primal and dual residuals in~\eqref{eq:convergence_check} can be rewritten in terms of the scaled problem as
\ifpreprint
\[
r_{\rm prim}^{k} = E^{-1}\bar{r}_{\rm prim}^{k} = E^{-1}(\bar{A}\bar{x}^{k} - \bar{z}^{k}), \qquad r_{\rm dual}^{k} = c^{-1}D^{-1}\bar{r}_{\rm dual}^{k} = c^{-1}D^{-1}(\bar{P}\bar{x}^{k} + \bar{q} + \bar{A}^{\tpose}\bar{y}^{k}),
\]
\else
\begin{align*}
r_{\rm prim}^{k} &= E^{-1}\bar{r}_{\rm prim}^{k} = E^{-1}(\bar{A}\bar{x}^{k} - \bar{z}^{k}), \\
r_{\rm dual}^{k} &= c^{-1}D^{-1}\bar{r}_{\rm dual}^{k} = c^{-1}D^{-1}(\bar{P}\bar{x}^{k} + \bar{q} + \bar{A}^{\tpose}\bar{y}^{k}),
\end{align*}
\fi
and the tolerances levels as
\begin{align*}
    \eps_{\rm prim} &= \eps_{\rm abs} + \eps_{\rm rel} \max\lbrace \|E^{-1}\bar{A} \bar{x}^{k}\|_{\infty}, \| E^{-1}\bar{z}^{k} \|_{\infty} \rbrace \\
    \eps_{\rm dual} &= \eps_{\rm abs} + \eps_{\rm rel} c^{-1}\max\lbrace \| D^{-1}\bar{P} \bar{x}^{k} \|_{\infty}, \| D^{-1}\bar{A}^\tpose \bar{y}^{k} \|_{\infty}, \| D^{-1}\bar{q} \|_{\infty} \rbrace.
\end{align*}

\paragraph{Quadratic programs infeasibility.}
When $\mathcal{C}=[l, u]$, the primal infeasibility conditions become
\begin{equation*}
	\left\|D^{-1}\bar{A}^\tpose \delta \bar{y}^{k} \right\|_{\infty} \le \eps_{\rm pinf} \| E \delta \bar{y}^{k}\|_{\infty},
		\quad
		\bar{u}^\tpose (\delta \bar{y}^{k})_{+} + \bar{l}^\tpose (\delta \bar{y}^{k})_{-} \le \eps_{\rm pinf} \| E \delta \bar{y}^{k}\|_{\infty},
\end{equation*}
where the primal infeasibility certificate is $c^{-1}E\delta \bar{y}^{k}$.
The dual infeasibility criteria are
\begin{gather*}
		\| D^{-1}\bar{P} \delta \bar{x}^{k} \|_{\infty} \le c\eps_{\rm dinf} \|D \delta \bar{x}^{k}\|_{\infty}, \quad
		\bar{q}^\tpose \delta \bar{x}^{k} \le c\eps_{\rm dinf}\|D \delta \bar{x}^{k}\|_{\infty}, \\[1em]
		(E^{-1}\bar{A} \delta \bar{x}^{k} )_i \begin{cases}
		\in \left[-\eps_{\rm dinf}, \eps_{\rm dinf}\right] \|D \delta \bar{x}^{k}\|_{\infty}& u_i, l_i \in \reals\\
		\ge -\eps_{\rm dinf}\|D \delta \bar{x}^{k}\|_{\infty} &u_i = +\infty\\
\le \eps_{\rm dinf}\|D \delta \bar{x}^{k}\|_{\infty} &l_i = -\infty,\end{cases}
\end{gather*}
where the dual infeasibility certificate is $D\delta \bar{x}^{k}$.

\subsection{Parameter selection}
\label{sec:parameter_selection}

The choice of parameters $(\rho,\sigma,\alpha)$ in Algorithm~\ref{alg:osqp_algorithm} is a key factor in determining the number of iterations required to find an optimal solution.
Unfortunately, it is still an open research question how to select the optimal \gls{ADMM} parameters, see~\cite{6892987,Nishihara:2015,7465685}.
After extensive numerical testing on millions of problem instances and a wide range of dimensions, we chose the algorithm parameters as follows for \glspl{QP}.

\paragraph{Choosing $\sigma$ and $\alpha$.}
The parameter~$\sigma$ is a regularization term which is used to ensure that a unique solution of~\eqref{eq:ADMM_iter1} will always exist, even when $P$ has one or more zero eigenvalues.
After scaling $P$ in order to minimize its condition number, we choose $\sigma$ as small as possible to preserve numerical stability without slowing down the algorithm.
We set the default value as $\sigma = 10^{-6}$.
The relaxation parameter~$\alpha$ in the range $[1.5, 1.8]$ has empirically shown to improve the convergence rate~\cite{Eckstein1994,doi:10.1287/ijoc.10.2.218}.
In the proposed method, we set the default value of~$\alpha=1.6$.

\paragraph{Choosing $\rho$.}

The most crucial parameter is the step-size~$\rho$.
Numerical testing showed that having different values of $\rho$ for different constraints, can greatly improve the performance.
For this reason, without altering the algorithm steps, we chose $\rho\in\symm_{++}^{m}$ being a positive definite diagonal matrix with different elements $\rho_i$.

For a specific \gls{QP}, if we know the active and inactive constraints, then we can rewrite it simply as an equality constrained \gls{QP}.
In this case the optimal $\rho$ is defined as $\rho_i=\infty$ for the active constraints and $\rho_i=0$ for the inactive constraints, therefore reducing the linear system~\eqref{eq:kkt} to the optimality conditions of the equivalent equality constrained \gls{QP} (after setting $\sigma = 0$).
Unfortunately, it is impossible to know a priori whether any given constraint is active or inactive at optimality, so we must instead adopt some heuristics.
We define $\rho$ as follows
\begin{equation*}
	\rho = \diag(\rho_1, \dots,\rho_m),\qquad \rho_i = \begin{dcases}\bar{\rho} & l_i\neq u_i\\
	10^3\bar{\rho} & l_i = u_i,\end{dcases}
\end{equation*}
where $\bar{\rho} > 0$.
In this way we assign a high value to the step-size related to the equality constraints since they will be active at the optimum.
Having a fixed value of $\bar{\rho}$ cannot provide fast convergence for different kind of problems since the optimal solution and the active constraints vary greatly.
To compensate for this issue, we adopt an adaptive scheme which updates $\bar{\rho}$ during the iterations based on the ratio between primal and dual residuals.
The idea of introducing ``feedback'' in the algorithm steps makes \gls{ADMM} more robust to bad scaling in the data; see~\cite{he2000,boyd2011distributed,wohlberg2017}.
Contrary to the adaptation approaches in the literature where the update increases or decreases the value of the step-size by a fixed amount, we adopt the following rule
\begin{equation*}
	\bar{\rho}^{k+1} \gets \bar{\rho}^{k} \sqrt{\frac{\|\bar{r}_{\rm prim}^{k}\|_{\infty} / \max\lbrace \|\bar{A}\bar{x}^{k}\|_{\infty}, \| \bar{z}^{k} \|_{\infty} \rbrace }{\|\bar{r}_{\rm dual}^{k}\|_{\infty} / \max\lbrace \| \bar{P} \bar{x}^{k} \|_{\infty}, \| \bar{A}^\tpose \bar{y}^{k} \|_{\infty}, \| \bar{q} \|_{\infty} \rbrace}}.
\end{equation*}
In other words we update $\bar{\rho}^{k}$ using the square root of the ratio between the scaled residuals normalized by the magnitudes of the relative part of the tolerances.
We set the initial value as $\bar{\rho}^{0}=0.1$.
In our benchmarks, if $\bar{\rho}^0$ does not already give a low number of \gls{ADMM} iterations, it gets usually tuned with a maximum of 1 or 2 updates.
The adaptation causes the KKT matrix in~\eqref{eq:kkt} to change and, if the linear system solver solution method is direct, it requires a new numerical factorization.
We do not require a new symbolic factorization because the sparsity pattern of the KKT matrix does not change.
Since the numerical factorization can be costly, we perform the adaptation only when it is really necessary.
In particular, we allow an update if the accumulated iterations time is greater than a certain percentage of the factorization time (nominally $40\:\%$) and if the new parameter is sufficiently different than the current one, \ie, $5$ times larger or smaller.
Note that in the case of an indirect method this rule allows for more frequent changes of $\rho$ since there is no need to factor the KKT matrix and the update is numerically much cheaper.
Note that the convergence of the \gls{ADMM} algorithm is hard to prove in general if the $\rho$ updates happen at each iteration.
However, if we assume that the updates stop after a fixed number of iterations the convergence results hold~\cite[Section 3.4.1]{boyd2011distributed}.

%% file: sections/parametric_programs.tex

\section{Parametric programs}
\label{sec:parametric_programs}
In application domains such as control, statistics, finance, and \gls{SQP}, problem~\eqref{pb:main} is  solved repeatedly for varying data.
For these problems, usually referred to as {\em parametric programs}, we can speed up the repeated OSQP calls by re-using the computations across multiple solves.

We make the distinction between cases in which only the vectors or all data in~\eqref{pb:main} change between subsequent problem instances.
We assume that the problem dimensions $n$ and $m$ and the sparsity patterns of $P$ and $A$ are fixed.

\paragraph{Vectors as parameters.}
If the vectors $q$, $l$, and $u$ are the only parameters that vary, then the KKT coefficient matrix in Algorithm~\ref{alg:osqp_algorithm} does not change across different instances of the parametric program.
Thus, if a direct method is used, we perform and store its factorization only once before the first solution and reuse it across all subsequent iterations.
Since the matrix factorization is the computationally most expensive step of the algorithm, this approach reduces significantly the amount of time OSQP takes to solve subsequent problems.
This class of problems arises very frequently in many applications including linear \gls{MPC} and \gls{MHE}~\cite{Rawlings:2009,Allgower1999}, Lasso~\cite{Tibshirani:1996,Candes:2008}, and portfolio optimization~\cite{OPT-001,JOFI:JOFI1525}.

\paragraph{Matrices and vectors as parameters.}
We separately consider the case in which the values (but not the locations) of the nonzero entries of matrices $P$ and $A$ are updated.
In this case, in a direct method, we need to refactor the matrix in Algorithm~\ref{alg:osqp_algorithm}.
However, since the sparsity pattern does not change we need only to recompute the {\em numerical factorization} while reusing the {\em symbolic factorization} from the previous solution.
This results in a modest reduction in the  computation time.
This class of problems encompasses several applications such as nonlinear \gls{MPC} and \gls{MHE}~\cite{Diehl2009} and sequential quadratic programming~\cite{NW06}.

\paragraph{Warm starting.}
In contrast to interior-point methods, OSQP is easily initialized by providing an initial guess of both the primal and dual solutions to the \gls{QP}.
This approach is known as \emph{warm starting} and is particularly effective when the subsequent \gls{QP} solutions do not vary significantly, which is the case for most {\em parametric programs} applications.
We can warm start the~\gls{ADMM} iterates from the previous OSQP solution $(x^\star, y^\star)$ by setting $(x^{0}, z^{0}, y^{0})\gets (x^{\star}, Ax^{\star}, y^{\star})$.
Note that we can warm-start the $\rho$ estimation described in Section~\ref{sec:osqp_implementation} to exploit the ratio between the primal and dual residuals to speed up convergence in subsequent solves.

%% file: sections/osqp_implementation.tex

\section{OSQP}
\label{sec:osqp_implementation}

We have implemented our proposed approach in the  ``Operator Splitting Quadratic Program'' (OSQP) solver, an open-source software package in the C language.
OSQP can solve any \gls{QP} of the form~\eqref{pb:qp} and makes no assumptions about the problem data other than convexity.
OSQP is available online at
\begin{center}
	\url{https://osqp.org}.
\end{center}
Users can call OSQP from C, \Cpp, Fortran, Python, Matlab, R, Julia, Ruby and Rust, and via parsers such as CVXPY~\cite{cvxpy,agrawal2018}, JuMP~\cite{DunningHuchetteLubin2017}, and YALMIP~\cite{Lofberg2004}.

To exploit the data sparsity pattern, OSQP accepts matrices in Compressed-Sparse-Column (CSC) format~\cite{davis2006direct}.
We implemented the linear system solution described in Section~\ref{sec:solving_linear_system} as an object-oriented interface to easily switch between efficient algorithms.
At present, OSQP ships with the open-source QDLDL direct solver which is our independent implementation based on \cite{Davis:2005}, and also supports dynamic loading of more advanced algorithms such as the MKL Pardiso direct solver~\cite{intel2017}. We plan to add iterative indirect solvers and other direct solvers in future versions.

The default values for the OSQP termination tolerances described in Section~\ref{sec:termination criteria} are
\[
\eps_{\rm abs} = \eps_{\rm rel} = 10^{-3},\quad \eps_{\rm pinf} = \eps_{\rm dinf} = 10^{-4}.
\]
The default step-size parameter $\sigma$ and the relaxation parameter $\alpha$ are set to
\[
	\sigma = 10^{-6},\quad \alpha = 1.6,
\]
while $\rho$ is automatically chosen by default as described in Section~\ref{sec:parameter_selection}, with optional user override.
We set the default fixed number of iterative refinement steps to 3.

OSQP reports the total computation time divided by the time required to perform preprocessing operations such as scaling or matrix factorization and the time to carry out the \gls{ADMM} iterations.
If the solver is called multiple times reusing the same matrix factorization, it will report only the \gls{ADMM} solve time as total computation time.
For more details we refer the reader to the solver documentation on the OSQP project website.

%% file: sections/examples.tex
\section{Numerical examples}
\label{sec:examples}
We benchmarked OSQP against the open-source interior-point solver ECOS~\cite{domahidi2013ecos}, the open-source active-set solver qpOASES~\cite{Ferreau:2014}, and the commercial interior-point solvers GUROBI~\cite{gurobi} and MOSEK~\cite{mosek}.
We executed every benchmark comparing different solvers with both low accuracy, \ie, $\eps_{\rm abs} = \eps_{\rm rel}  = 10^{-3}$, and high accuracy, \ie, $\eps_{\rm abs} = \eps_{\rm rel}  = 10^{-5}$.
We set GUROBI, ECOS, MOSEK and OSQP  primal and dual feasibility tolerances to our low and high accuracy tolerances. Since qpOASES is an active-set method and does not allow the user to tune primal nor dual feasibility tolerances, we set it to its default termination settings.
In addition, the maximum time we allow each solver to run is $1000\; {\rm sec}$ and no limit on the maximum number of iterations.
Note that the use of maximum time limits with no bounds on the number of iterations is the default setting in commercial solvers such as MOSEK.
For every solver we leave all the other settings to the internal defaults.

In general it is hard to compare the solution accuracies because all the solvers, especially commercial ones, use an internal problem scaling and verify that the termination conditions are satisfied against their scaled version of the problem. In contrast, OSQP allows the option to check the termination conditions against the internally scaled or the original problem. Therefore, to make the benchmark fair, we say that the primal-dual solution $(x^\star, y^\star)$ returned by each solver is optimal if the following optimality conditions are satisfied with tolerances defined above with low and high accuracy modes,
\begin{align*}
	& \|(Ax^\star - u)_{+} + (Ax^\star - l)_{-}\|_\infty \leq \eps_{\rm prim},
	&&\|Px^\star + q + A^\tpose y^\star\|_\infty \leq \eps_{\rm dual},
\end{align*}
where $\eps_{\rm prim}$ and $\eps_{\rm dual}$ are defined in Section~\ref{sec:termination criteria}.
If the primal-dual solution returned by a solver does not satisfy the optimality conditions defined above, we consider it a failure.
Note that we decided not to include checks on the complementary slackness satisfaction because interior-point solvers satisfied them with different metrics and scalings, therefore failing very often. In contrast OSQP always satisfies complementary slackness conditions with machine precision by construction.

In addition, we used the direct single-threaded linear system solver QDLDL~\cite{qdldl} based on~\cite{Amestoy:2004,Davis:2005} and very simple linear algebra where other solvers such as GUROBI and MOSEK use advanced multi-threaded linear system solvers and custom linear algebra.

All the experiments were carried out on the MIT SuperCloud facility in collaboration with the Lincoln Laboratory~\cite{supercloud} with 16 Intel Xeon E5-2650 cores.
The code for all the numerical examples is available online at~\cite{OSQPBenchmarks}.

\paragraph{Shifted geometric mean.}
As in most common benchmarks~\cite{hansbench}, we make use of the normalized shifted geometric mean to compare the timings of the various solvers.
Given the time required by solver $s$ to solve problem $p$ $t_{p, s}$, we define the shifted geometric mean as
\begin{equation*}
g_{s} \eqdef \sqrt[n]{\prod_{p} (t_{p, s} + k) - k},
\end{equation*}
where $n$ is the number of problem instances considered and $k=1$ is the shift~\cite{hansbench}.
The normalized shifted geometric mean is therefore
\begin{equation*}
r_{s} \eqdef g_{s} / \min_{s} g_{s}.
\end{equation*}
This value shows the factor at which a specific solver is slower than the fastest one with scaled value of $1.00$.
If solver $s$ fails at solving problem $p$, we set the time as the maximum allowed, \ie, $t_{p, s} = 1000\; {\rm sec}$.
Note that to avoid memory overflows in the product, we compute in practice the shifted geometric mean as $e^{\ln{g_s}}$.

\paragraph{Performance profiles.}
We also make use of the performance profiles~\cite{Dolan2002} to compare the solver timings.
We define the performance ratio
\begin{equation*}
	u_{p, s} \eqdef t_{p, s} /\min_{s} t_{p, s}.
\end{equation*}
The performance profile plots the function $f_s: \reals \mapsto [0, 1]$ defined as
\begin{equation*}
	f_{s}(\tau) \eqdef \frac{1}{n} \sum_{p}\mathcal{I}_{\leq \tau}(u_{p, s}), 
\end{equation*}
where $\mathcal{I}_{\leq \tau}(u_{p, s}) = 1$ if $u_{p, s} \leq \tau$ or $0$ otherwise.
The value $f_s(\tau)$ corresponds to the fraction of problems solved within $\tau$ times from the best solver.
Note that while we cannot necessarily assess the performance of one solver relative to another with performance profiles, they still represent a viable choice to benchmark the performance of a solver with respect to the best one~\cite{gould2016}.

\subsection{Benchmark problems}
We considered \glspl{QP} in the form~\eqref{pb:qp} from 7 problem classes ranging from standard random programs to applications in the areas of control, portfolio optimization and machine learning.
For each problem class, we generated $10$ different instances for $20$ dimensions giving a total of $1400$ problem instances.
All instances were obtained from either real data or from non-trivial random data.
Note that the random \glspl{QP} and random equality constrained \glspl{QP} problem classes might not closely correspond to a real-world application.
However, they have a typical number of nonzero elements appearing in practice.
We described generation for each class in Appendix~\ref{app:problem_classes}.
Throughout all the problem classes, $n$ ranges between $10^1$ and $10^4$, $m$ between $10^2$ and $10^5$, and the number of nonzeros $N$ between $10^2$ and $10^8$.

\paragraph{Results.}
We show in Figures~\ref{fig:computation_times_low_accuracy} and~\ref{fig:computation_times_high_accuracy} the OSQP and GUROBI computation times across all the problem classes for low and high accuracy solutions respectively.
OSQP is competitive or even faster than GUROBI for several problem classes.
Results are shown in Table~\ref{tab:benchmark_results} and Figure~\ref{fig:benchmarks_performance_profiles}.
OSQP shows the best performance across these benchmarks with MOSEK performing better at lower accuracy and GUROBI at higher accuracy. ECOS is generally slower than the other interior-point solvers but faster than qpOASES that shows issues with many constraints.
Table~\ref{tab:benchmark_osqp_stats} contains the OSQP statistics for this benchmark class.
Because of the good convergence behavior of OSQP on these problems, the setup time is significant compared to the solve time, especially at low accuracy.
Solution polishing increases the solution time by a median of $10$ to $20$ percent due to the additional factorization used.
The worst-case time increase is very high and happens for the problems that converge in very few iterations.
Note that with high accuracy, polishing succeeds in 83\% of test cases while on low accuracy it succeeds in only 42\% of cases.
The number of $\rho$ updates is in general very low, usually requiring just more matrix factorization to adjust, with up to 5 refactorisations used in the worst case when solving with high accuracy.
\begin{table}
    \caption{Benchmark problems comparison with timings as shifted geometric mean and failure rates.}
\label{tab:benchmark_results}
      \centering
      \tablecompareosqpvsall{benchmark_problems}
\end{table}

\begin{table}
    \caption{Benchmark problems OSQP statistics.}
\label{tab:benchmark_osqp_stats}
      \centering
      \osqpstatstable{benchmark_problems}
\end{table}

\subsection{SuiteSparse matrix collection least squares problems}

We considered 30 least squares problem in the form $Ax \approx b$ from the SuiteSparse Matrix Collection library~\cite{davis2011}.
Using the Lasso and Huber problem setups from Appendix~\ref{app:problem_classes} we formulate 60 \glspl{QP} that we solve with OSQP, GUROBI and MOSEK.
We excluded ECOS because its interior-point algorithm showed numerical issues for several problems of the test set.
We also excluded qpOASES because it is not designed for large linear systems.

\paragraph{Results.}
Results are shown in Table~\ref{tab:suitesparse_results} and Figure~\ref{fig:suitesparse_performance_profiles}.
OSQP shows the best performance with GUROBI slightly slower and MOSEK third.
The failure rates for GUROBI and MOSEK are higher because the reported solution does not satisfy the optimality conditions of the original problem.
We display the OSQP statistics in Table~\ref{tab:suitesparse_osqp_stats}.
The setup phase takes a significant amount of time compared to the solve phase, especially when OSQP converges in a few iterations.
This happens because the large problem dimensions result in a large initial factorization time.
Polish time is in general $22$ to $32\%$ of the total solution time.
However, the success is usually reliable, succeeding $78\%$ of the times with very high quality solutions.
The number of matrix refactorizations required due to $\rho$ updates is very low in these examples, with a maximum of $2$ or $3$ even for high accuracy.

\begin{table}
	\caption{SuiteSparse matrix problems comparison with timings as shifted geometric mean and failure rates.}
	\label{tab:suitesparse_results}
	\centering
	\tablecompareosqpvsip{suitesparse_problems}
\end{table}

\begin{table}
	\caption{SuiteSparse problems OSQP statistics.}
	\label{tab:suitesparse_osqp_stats}
	\centering
	\osqpstatstable{suitesparse_problems}
\end{table}

\subsection{Maros-M\'{e}sz\'{a}ros problems}
We considered the Maros-M\'{e}sz\'{a}ros test set~\cite{maros1999} of hard \glspl{QP}.
We compared the OSQP solver against GUROBI and MOSEK against all the problems in the set.
We decided to exclude ECOS because its interior-point algorithm showed numerical issues for several problems of the test set.
We also excluded qpOASES because it could not solve most of the problems since it is not suited for large \glspl{QP} -- it is based on an active-set method with dense linear algebra.

\paragraph{Results.}
Results are shown in Table~\ref{tab:maros_results} and Figure~\ref{fig:maros_performance_profiles}.
GUROBI shows the best performance and OSQP, while slower, is still competitive on both low and high accuracy tests.
MOSEK remains the slowest in every case. Table~\ref{tab:maros_osqp_stats} shows the statistics relative to OSQP.
Since these hard problems require a larger number of iterations to converge, the setup time overhead compared to the solution time is in general lower than the other benchmark sets.
Moreover, since the problems are badly scaled and degenerate, the polishing strategy rarely succeeds.
However, the median time increase from the polish step is less than $10\%$ of the total computation time for both low and high accuracy modes.
Note that the number of $\rho$ updates is usually very low with a median of $1$ or $2$. However, there are some worst-case problems when it is very high because the bad scaling causes issues in our $\rho$ estimation. However, from our data we have seen that in more than $95\%$ of the cases the number of $\rho$ updates is less than $5$.

\begin{table}
	\caption{Maros-M\'{e}sz\'{a}ros problems comparison with timings as shifted geometric mean and failure rates.}
	\label{tab:maros_results}
	\centering
	\tablecompareosqpvsip{maros_meszaros_problems}
\end{table}

\begin{table}
	\caption{Maros-M\'{e}sz\'{a}ros problems OSQP statistics.}
	\label{tab:maros_osqp_stats}
	\centering
	\osqpstatstable{maros_meszaros_problems}
\end{table}

\subsection{Warm start and factorization caching}
To show the benefits of warm starting and factorization caching, we solved a sequence of \glspl{QP} using OSQP with the data varying according to some parameters. Since we are not comparing OSQP with other high accuracy solvers in these benchmarks, we use its default settings with accuracy $10^{-3}$.

\paragraph{Lasso regularization path.}
We solved a Lasso problem described in Appendix~\ref{app:lasso} with varying $\lambda$ in order to choose a regressor with good validation set performance.
We solved one problem instance with $n=50, 100, 150, 200$ features, $m=100n$ data points, and $\lambda$ logarithmically spaced taking 100 values between $\lambda_{\rm max} = \|A^\tpose b\|_{\infty}$ and $0.01\lambda_{\rm max}$.

Since the parameters only enter linearly in the cost, we can reuse the matrix factorization and enable warm starting to reduce the computation time as discussed in Section~\ref{sec:parametric_programs}.

\paragraph{Model predictive control.}
In \gls{MPC}, we solve the optimal control problem described in Appendix~\ref{app:optimal_control} at each time step to compute an optimal input sequence over the horizon.
Then, we apply only the first input to the system and propagate the state to the next time step.
The whole procedure is repeated with an updated initial state $x_{\rm init}$.
We solved the control problem with $n_x = 20, 40, 60, 80$ states, $n_u=n_x/2$ inputs, horizon $T=10$ and $100$ simulation steps.
The initial state of the simulation is uniformly distributed and constrained to be within the feasible region, \ie, $x_{{\rm init}} \sim \mathcal{U}(-0.5\overline{x}, 0.5\overline{x})$.

Since the parameters only enter linearly in the constraints bounds, we can reuse the matrix factorization and enable warm starting to reduce the computation time as discussed in Section~\ref{sec:parametric_programs}.

\paragraph{Portfolio back test.}
Consider the portfolio optimization problem in Appendix~\ref{app:portfolio} with $n=10k$ assets and $k=100, 200, 300, 400$ factors.

We run a 4 years back test to compute the optimal assets investment depending on varying expected returns and factor models~\cite{OPT-023}.
We solved $240$ \glspl{QP} per year giving a total of 960 \glspl{QP}.
Each month we solved $20$ \glspl{QP} corresponding to the trading days.
Every day, we updated the expected returns $\mu$ by randomly generating another vector with $\mu_i \sim 0.9 \hat{\mu}_i + \mathcal{N}(0,0.1)$, where $\hat{\mu}_i$ comes from the previous expected returns.
The risk model was updated every month by updating the nonzero elements of $D$ and $F$ according to $D_{ii} \sim 0.9 \hat{D}_{ii} + \mathcal{U}[0, 0.1\sqrt{k}]$ and $F_{ij} \sim 0.9\hat{F}_{ij} + \mathcal{N}(0,0.1)$ where $\hat{D}_{ii}$ and $\hat{F}_{ij}$ come from the previous risk model.

As discussed in Section~\ref{sec:parametric_programs}, we exploited the following computations during the \gls{QP} updates to reduce the computation times.
Since $\mu$ only enters in the linear part of the objective, we can reuse the matrix factorization and enable warm starting.
Since the sparsity patterns of $D$ and $F$ do not change during the monthly updates, we can reuse the symbolic factorization and exploit warm starting to speed up the computations.

\paragraph{Results.}
We show the results in Table~\ref{tab:parametric}.
For the Lasso problem we see more than $10$-fold improvement in time and between $8$ and $11$ times reduction in number of iterations depending on the dimension.
For the MPC problem the number of iterations does not significantly decrease because the number of iterations is already low in cold-start. However we get from $2.6$ to $4$-fold time improvement from factorization caching.
OSQP shows from $5.8$ to $7$ times reduction in time for the portfolio problem and from $2.9$ to $3.6$ times reduction in number of iterations.

\begin{table}[]
	\centering
	\caption{OSQP parametric problem results with warm start (ws) and without warm start (no ws) in terms of time in seconds and number of iterations for different leading problem dimensions of Lasso, MPC and Portfolio classes.}
	\label{tab:parametric}
	\begin{tabular}{@{}rr
			S[table-figures-integer=1]
			S[table-figures-integer=1]
			S[table-figures-integer=2]
			S[table-figures-integer=3]
			S[table-figures-integer=3]
			S[table-figures-integer=2]
			@{}
		}
		\toprule
		Problem & dim. & $\shortstack{Time\\no ws}$ & $\shortstack{Time\\ws}$ & $\shortstack{Time\\improv.}$ & $\shortstack{Iter\\no ws}$ & $\shortstack{Iter\\ws}$ & $\shortstack{Iter\\improv.}$\\
		\midrule
		\multirow{4}*{Lasso}
	\begin{filecontents*}{./data/results/parametric_problems/statistics_lasso.csv}
,problem,dimension,OSQP_ws_mean_time,OSQP_ws_median_time,OSQP_ws_mean_iter,OSQP_ws_median_iter,OSQP_nows_mean_time,OSQP_nows_median_time,OSQP_nows_mean_iter,OSQP_nows_median_iter,time_speedup_mean,time_speedup_median,iter_reduction_mean,iter_reduction_median
0,Lasso,50,0.011626220609999998,0.009811594,25.75,25.0,0.22500497389,0.22101784200000002,210.25,200.0,19.35323450653153,22.526191157114738,8.16504854368932,8.0
1,Lasso,100,0.04011979855,0.03824527,25.75,25.0,0.42349856996,0.4277865835,224.0,225.0,10.555849861314918,11.185346148687145,8.699029126213592,9.0
2,Lasso,150,0.08597965517000002,0.081400697,25.75,25.0,1.02197516542,1.0998814344999999,235.5,275.0,11.886244058543133,13.511941236817664,9.145631067961165,11.0
3,Lasso,200,0.14935050164999997,0.140244054,26.0,25.0,2.08881798686,2.1762049185,281.75,300.0,13.986012526125322,15.517270475509784,10.836538461538462,12.0
\end{filecontents*}
\csvreader[head to column names, late after line=\\]{./data/results/parametric_problems/statistics_lasso.csv}{
		dimension=\n,
		OSQP_ws_mean_time=\timews,
		OSQP_nows_mean_time=\timenows,
		time_speedup_mean=\timespeedup,
		OSQP_ws_mean_iter=\iterws,
		OSQP_nows_mean_iter=\iternows,
		iter_reduction_mean=\iterreduct,
	}
	{& \n & \timenows & \timews & \timespeedup & \iternows & \iterws & \iterreduct}
	\midrule
	\multirow{4}*{MPC}
\begin{filecontents*}{./data/results/parametric_problems/statistics_mpc.csv}
,problem,dimension,OSQP_ws_mean_time,OSQP_ws_median_time,OSQP_ws_mean_iter,OSQP_ws_median_iter,OSQP_nows_mean_time,OSQP_nows_median_time,OSQP_nows_mean_iter,OSQP_nows_median_iter,time_speedup_mean,time_speedup_median,iter_reduction_mean,iter_reduction_median
0,MPC,20,0.0017661863900000003,0.001329298,32.75,25.0,0.0071022403300000005,0.0074996995,89.5,100.0,4.021229225982202,5.641849683065798,2.732824427480916,4.0
1,MPC,40,0.00507617565,0.0045008245,27.25,25.0,0.013662178799999997,0.012491060499999998,29.0,25.0,2.6914314519435507,2.7752827287533646,1.0642201834862386,1.0
2,MPC,60,0.01295097151,0.009368037,33.0,25.0,0.03461645991,0.0304179585,33.75,25.0,2.6728851872827573,3.246993847270245,1.0227272727272727,1.0
3,MPC,80,0.021680080979999995,0.016191612,31.75,25.0,0.0667613624,0.06008841299999999,32.0,25.0,3.079387132436809,3.7110828125081055,1.0078740157480315,1.0
\end{filecontents*}
\csvreader[head to column names, late after line=\\]{./data/results/parametric_problems/statistics_mpc.csv}{
	dimension=\n,
	OSQP_ws_mean_time=\timews,
	OSQP_nows_mean_time=\timenows,
	time_speedup_mean=\timespeedup,
	OSQP_ws_mean_iter=\iterws,
	OSQP_nows_mean_iter=\iternows,
	iter_reduction_mean=\iterreduct,
}
{& \n & \timenows & \timews & \timespeedup & \iternows & \iterws & \iterreduct}
\midrule
\multirow{4}*{Portfolio}
\begin{filecontents*}{./data/results/parametric_problems/statistics_portfolio.csv}
,problem,dimension,OSQP_ws_mean_time,OSQP_ws_median_time,OSQP_ws_mean_iter,OSQP_ws_median_iter,OSQP_nows_mean_time,OSQP_nows_median_time,OSQP_nows_mean_iter,OSQP_nows_median_iter,time_speedup_mean,time_speedup_median,iter_reduction_mean,iter_reduction_median
0,Portfolio,100,0.030487230036458333,0.028728840999999998,25.416666666666668,25.0,0.17735215724479167,0.172825838,93.33333333333333,100.0,5.817260440935567,6.015760886420724,3.6721311475409832,4.0
1,Portfolio,200,0.06056596586041667,0.054893913499999995,25.390625,25.0,0.41614651463541663,0.406007738,86.875,75.0,6.870963068507626,7.396225047791502,3.4215384615384616,3.0
2,Portfolio,300,0.09739454775833334,0.0877059675,25.520833333333332,25.0,0.6462597646145833,0.6226212615,80.52083333333333,75.0,6.635481959607822,7.09896121378514,3.1551020408163266,3.0
3,Portfolio,400,0.13937640847083335,0.125572368,26.09375,25.0,0.976060201328125,0.91235479,76.45833333333333,75.0,7.003051750557775,7.265569683292108,2.930139720558882,3.0
\end{filecontents*}
\csvreader[head to column names, late after line=\\]{./data/results/parametric_problems/statistics_portfolio.csv}{
	dimension=\n,
	OSQP_ws_mean_time=\timews,
	OSQP_nows_mean_time=\timenows,
	time_speedup_mean=\timespeedup,
	OSQP_ws_mean_iter=\iterws,
	OSQP_nows_mean_iter=\iternows,
	iter_reduction_mean=\iterreduct,
}
{& \n & \timenows & \timews & \timespeedup & \iternows & \iterws & \iterreduct}
\bottomrule
\end{tabular}
\end{table}

%% file: sections/conclusion.tex

\section{Conclusions}
\label{sec:conclusions}

We presented a novel general-purpose \gls{QP} solver based on \gls{ADMM}.
Our method uses a new splitting requiring the solution of a quasi-definite linear system that is always solvable independently from the problem data.
We impose no assumptions on the problem data other than convexity, resulting in a general-purpose and very robust algorithm.

For the first time, we propose a first-order \gls{QP} solution method able to provide primal and dual infeasibility certificates if the problem is unsolvable without resorting to homogeneous self-dual embedding or additional complexity in the iterations.

In contrast to other first-order methods, our solver can provide high-quality solutions by performing \emph{solution polishing}.
After guessing which constraints are active, we compute the solutions of an additional small equality constrained \gls{QP} by solving a linear system.
If the constraints are identified correctly, the returned solution has accuracy equal or higher than interior-point methods.

The proposed method is easily warm started to reduce the number of iterations.
If the problem matrices do not change, the linear system matrix factorization can be cached and reused across multiple solves greatly improving the computation time.
This technique can be extremely effective, especially when solving parametric \glspl{QP} where only part of the problem data change.

We have implemented our algorithm in the open-source OSQP solver written in C and interfaced with multiple other languages and parsers.
OSQP is based on sparse linear algebra and is able to exploit the structure of \glspl{QP} arising in different application areas.
OSQP is robust against noisy and unreliable data and, after the first factorization is computed, can be compiled to be library-free and division-free, making it suitable for embedded applications.
Thanks to its simple and parallelizable iterations, OSQP can handle large-scale problems with millions of nonzeros.

We extensively benchmarked the OSQP solver with problems arising in several application domains including finance, control and machine learning.
In addition, we benchmarked it against the hard problems from the Maros-M\'{e}sz\'{a}ros test set~\cite{maros1999} and Lasso and Huber fitting problems generated with sparse matrices from the SuiteSparse Matrix Collection~\cite{davis2011}.
Timing and failure rate results showed great improvements over state-of-the-art academic and commercial \gls{QP} solvers.

OSQP has already a large userbase with tens of thousands of users both from top academic institutions and large corporations.

%% file: sections/appendix.tex

\appendix
\section{Problem classes}\label{app:problem_classes}
In this section we describe the random problem classes used in the benchmarks and derive formulations with explicit linear equalities and inequalities that can be directly written in the form $Ax \in \mathcal{C}$ with $\mathcal{C} = [l, u]$.

\subsection{Random QP}
Consider the following \gls{QP}
\begin{equation*}
\begin{array}{ll}
\mbox{minimize}   & (1/2) x^\tpose P x + q^\tpose x \\
\mbox{subject to} & l \leq A x \leq u.
\end{array}
\end{equation*}

\paragraph{Problem instances.}
The number of variables and constraints in our problem instances are $n$ and $m=10n$.
We generated random matrix $P= M M^\tpose + \alpha I$ where $M\in \reals^{n \times n}$ and $15 \%$ nonzero elements $M_{ij}\sim \mathcal{N}(0,1)$.
We add the regularization $\alpha I$ with $\alpha = 10^{-2}$ to ensure that the problem is not unbounded.
We set the elements of $A\in \reals^{m \times n}$ as $A_{ij} \sim \mathcal{N}(0,1)$ with only $15 \%$ being nonzero.
The linear part of the cost is normally distributed, \ie, $q_i \sim  \mathcal{N}(0,1)$. We generated the constraint bounds as $u_i \sim \mathcal{U}(0,1)$, $l_i \sim -\mathcal{U}(0,1)$.

\subsection{Equality constrained QP} 
Consider the following equality constrained \gls{QP}
\begin{equation*}
\begin{array}{ll}
\mbox{minimize}   & (1/2) x^\tpose P x + q^\tpose x \\
\mbox{subject to} & A x = b.
\end{array}
\end{equation*}
This problem can be rewritten as~\eqref{pb:main} by setting $l = u = b$.

\paragraph{Problem instances.}
The number of variables and constraints in our problem instances are $n$ and $m=\lfloor n/2\rfloor$.

We generated random matrix $P= M M^\tpose + \alpha I$ where $M\in \reals^{n \times n}$ and $15 \%$ nonzero elements $M_{ij}\sim \mathcal{N}(0,1)$.
We add the regularization $\alpha I$ with $\alpha = 10^{-2}$ to ensure that the problem is not unbounded.
We set the elements of $A\in \reals^{m \times n}$ as $A_{ij} \sim \mathcal{N}(0,1)$ with only $15 \%$ being nonzero.
The vectors are all normally distributed, \ie, $q_i, b_i \sim  \mathcal{N}(0,1)$.

\paragraph{Iterative refinement interpretation.}
Solution of the above problem can be found directly by solving the following linear system
\begin{equation}\label{eq:eq_constr_qp_optimality}
    \begin{bmatrix} P & A^\tpose \\ A & 0\end{bmatrix} \begin{bmatrix}x \\ \nu \end{bmatrix} = \begin{bmatrix} -q \\ b \end{bmatrix}.
\end{equation}
If we apply the \gls{ADMM} iterations~\eqref{eq:ADMM_iter1}--\eqref{eq:ADMM_iter5} for solving the above problem, and by setting $\alpha=1$ and $y^0=b$, the algorithm boils down to the following iteration
\begin{equation*}
    \begin{bmatrix} x^{k+1} \\ \nu^{k+1} \end{bmatrix} = \begin{bmatrix} x^k \\ \nu^k \end{bmatrix} + \begin{bmatrix} P + \sigma I & A^\tpose \\ A & -\rho^{-1} I \end{bmatrix}^{-1}\left( \begin{bmatrix} -q \\ b \end{bmatrix} - \begin{bmatrix} P & A^\tpose \\ A & 0\end{bmatrix} \begin{bmatrix} x^k \\ \nu^k \end{bmatrix} \right),
\end{equation*}
which is equivalent to \eqref{eq:iterative_refinement} with $g = (-q, b)$ and $\hat{t}^k = (x^k,\nu^k)$.
This means that Algorithm~\ref{alg:osqp_algorithm} applied to solve an equality constrained \gls{QP} is equivalent to applying iterative refinement~\cite{Wilkinson:1963,duff1986direct} to solve the KKT system \eqref{eq:eq_constr_qp_optimality}.
Note that the perturbation matrix in this case is
\[
	\Delta K = \begin{bmatrix} \sigma I \\ & -\rho^{-1} I \end{bmatrix},
\]
which justifies using a low value of $\sigma$ and a high value of $\rho$ for equality constraints.

\subsection{Optimal control}\label{app:optimal_control}

We consider the problem of controlling a constrained linear time-invariant dynamical system.
To achieve this, we formulate the following optimization problem~\cite{borrelli2017predictive}
\begin{equation}\label{pb:mpc}
\begin{array}{ll}
	\mbox{minimize}  & x_T^\tpose Q_T x_T + \sum_{t=0}^{T-1} x_t^\tpose Q x_t + u_t^\tpose R u_t \\[.25em]
  \mbox{subject to}
  & x_{t+1} = Ax_t + Bu_t \\
  & x_t \in \mathcal{X}, u_t \in \mathcal{U} \\
  & x_0 = x_{\rm init}.
\end{array}
\end{equation}
The states $x_t\in\reals^{n_x}$ and the inputs $u_k\in\reals^{n_u}$ are subject to polyhedral constraints defined by the sets~$\mathcal{X}$ and~$\mathcal{U}$.
The horizon length is~$T$ and the initial state is $x_{\rm init}\in\reals^{n_x}$.
Matrices $Q\in\symm^{n_x}_+$ and $R\in\symm^{n_u}_{++}$  define the state and input costs at each stage of the horizon, and $Q_T\in\symm^{n_x}_+$ defines the final stage cost.

By defining the new variable $z = (x_0, \dots, x_{T}, u_0, \dots, u_{T-1})$, problem~\eqref{pb:mpc} can be written as a sparse QP of the form~\eqref{pb:qp} with a total of $n_x(T+1) + n_u T$ variables.

\paragraph{Problem instances.}
We defined the linear systems with $n = n_x$ states and $n_u = 0.5n_x$ inputs.
We set the horizon length to $T=10$.
We generated the dynamics as $A = I + \Delta$ with $\Delta_{ij} \sim \mathcal{N}(0, 0.01)$.
We chose only stable dynamics by enforcing the norm of the eigenvalues of $A$ to be less than $1$.
The input action is modeled as $B$ with $B_{ij} \sim \mathcal{N}(0, 1)$.

The state cost is defined as $Q = \diag(q)$ where $q_i \sim \mathcal{U}(0, 10)$ and $70\%$ nonzero elements in $q$.
We chose the input cost as $R = 0.1I$.
The terminal cost $Q_T$ is chosen as the optimal cost for the linear quadratic regulator (LQR) applied to $A, B, Q, R$ by solving a discrete algebraic Riccati equation (DARE)~\cite{borrelli2017predictive}.
We generated input and state constraints as
\begin{equation*}
    \mathcal{X} = \{x_t \in \reals^{n_x} \mid -\overline{x} \le x_t \le \overline{x} \}, \qquad \mathcal{U}=\{u_t \in \reals^{n_u} \mid -\overline{u} \le u_t \le \overline{u}\},
\end{equation*}
where  $\overline{x}_i\sim \mathcal{U}(1, 2)$ and $\overline{u}_i \sim \mathcal{U}(0, 0.1)$.
The initial state is uniformly distributed with $x_{{\rm init}} \sim \mathcal{U}(-0.5\overline{x}, 0.5\overline{x})$.

\subsection{Portfolio optimization}\label{app:portfolio}

Portfolio optimization is a problem arising in finance that seeks to allocate assets in a way that maximizes the risk adjusted return~\cite{OPT-001,JOFI:JOFI1525,OPT-023},~\cite[\S 4.4.1]{boyd2004convex},

\begin{equation*}
\begin{array}{ll}
	\mbox{maximize}   & \mu^\tpose x - \gamma (x^\tpose \Sigma x) \\
	\mbox{subject to} & \ones^\tpose x = 1 \\
	& x \ge 0,
\end{array}
\end{equation*}
where the variable $x \in \reals^{n}$ represents the portfolio, $\mu  \in \reals^{n}$ the vector of expected returns, $\gamma > 0$ the risk aversion parameter, and $\Sigma\in\symm_+^n$ the risk model covariance matrix.
The risk model is usually assumed to be the sum of a diagonal and a rank $k < n$ matrix
\begin{equation*}
    \Sigma = FF^\tpose + D,
\end{equation*}
where $F\in \reals^{n\times k}$ is the factor loading matrix and $D\in \reals^{n\times n}$ is a diagonal matrix describing the asset-specific risk.

We introduce a new variable $y=F^\tpose x$ and solve the resulting problem in variables~$x$ and~$y$
\begin{equation}\label{pb:portfolio_sparse}
\begin{array}{ll}
	\mbox{minimize}   & x^\tpose D x + y^\tpose y - \gamma^{-1} \mu^\tpose x \\
  \mbox{subject to}
  & y=F^\tpose x \\
  & \ones^\tpose x = 1 \\
	& x \ge 0,
\end{array}
\end{equation}
Note that the Hessian of the objective in \eqref{pb:portfolio_sparse} is a diagonal matrix.
Also, observe that $FF^\tpose$ does not appear in problem \eqref{pb:portfolio_sparse}.

\paragraph{Problem instances.} We generated portfolio problems for increasing number of factors $k$ and number of assets $n=100k$.
The elements of matrix $F$ were chosen as $F_{ij} \sim \mathcal{N}(0,1)$ with $50 \%$ nonzero elements.
The diagonal matrix $D$ is chosen as $D_{ii} \sim \mathcal{U}[0, \sqrt{k}]$.
The mean return was generated as $\mu_i \sim \mathcal{N}(0,1)$.
We set $\gamma = 1$.

\subsection{Lasso}\label{app:lasso}
The {\em least absolute shrinkage and selection operator (Lasso)} is a well known linear regression technique obtained by adding an~$\ell_1$ regularization term in the objective~\cite{Tibshirani:1996,Candes:2008}.
It can be formulated as
\begin{equation*}
\begin{array}{ll}
	\mbox{minimize}   & \|Ax - b\|_2^2 + \lambda\|x\|_1,
\end{array}
\end{equation*}
where $x\in \reals^{n}$ is the vector of parameters and $A\in\reals^{m \times n}$ is the data matrix and $\lambda$ is the weighting parameter.

We convert this problem to the following \gls{QP}
\begin{equation*}
\begin{array}{ll}
	\mbox{minimize}   & y^\tpose y + \lambda \ones^\tpose t\\
    \mbox{subject to} & y = Ax - b\\
    & -t \leq x \leq t,
\end{array}
\end{equation*}
where $y\in\reals^{m}$ and $t\in\reals^{n}$ are two newly introduced variables.

\paragraph{Problem instances.}
The elements of matrix $A$ are generated as $A_{ij} \sim \mathcal{N}(0,1)$ with $15 \%$ nonzero elements. To construct the vector $b$, we generated the true sparse vector $v\in \reals^{n}$ to be learned
\[
v_i \sim
\begin{dcases}
    0 & \mbox{with probability } p=0.5\\
    \mathcal{N}(0, 1/n) & \mbox{otherwise}.
\end{dcases}
\]
Then we let $b=Av + \varepsilon$ where $\varepsilon$ is the noise generated as $\varepsilon_i \sim \mathcal{N}(0, 1)$. We generated the instances with varying $n$ features and $m = 100n$ data points.
The parameter $\lambda$ is chosen as $(1/5)\|A^\tpose b\|_{\infty}$ since $\|A^\tpose b\|_{\infty}$ is the critical value above which the solution of the problem is $x=0$.

\subsection{Huber fitting}

\emph{Huber fitting} or the \emph{robust least-squares problem} performs linear regression under the assumption that there are outliers in the data~\cite{Huber:1964,Huber:1981}.
The fitting problem is written as
\begin{equation}\label{pb:huber_fitting}
\begin{array}{ll}
	\mbox{minimize}   & \sum_{i=1}^m \phi_{\rm hub}(a_i^T x - b_i),
\end{array}
\end{equation}
with the Huber penalty function $\phi_{\rm hub}:\reals\to\reals$ defined as
\begin{equation*}
    \phi_{\rm hub}(u) = \begin{cases}
        u^2         & |u| \le M \\
        M(2|u|-M)   & |u| > M.
    \end{cases}
\end{equation*}
Problem~\eqref{pb:huber_fitting} is equivalent to the following \gls{QP}~\cite[Eq.~(24)]{Mangasarian:2000}
\begin{equation*}
	\begin{array}{ll}
		\mbox{minimize}     & u^\tpose u + 2M \ones^\tpose (r + s) \\
		\mbox{subject to}   & Ax - b - u = r - s \\
				    & r \ge 0 \\
				    & s \ge 0.
	\end{array}
\end{equation*}

\paragraph{Problem instances.}
We generate the elements of $A$ as $A_{ij} \sim \mathcal{N}(0,1)$ with $15 \%$ nonzero elements.
To construct $b\in\reals^m$ we first generate a vector $v\in\reals^n$ as $v_i \sim \mathcal{N}(0,1/n)$ and a noise vector $\varepsilon\in\reals^m$ with elements
\begin{equation*}
    \varepsilon_i \sim \begin{cases}
        \mathcal{N}(0,1/4)  & \text{with probability $p=0.95$} \\
        \mathcal{U}[0,10]   & \text{otherwise}.
    \end{cases}
\end{equation*}
We then set $b = Av + \varepsilon$.
For each instance we choose $m=100n$ and $M=1$.

\subsection{Support vector machine}

\emph{Support vector machine} problem seeks an affine function that approximately classifies the two sets of points~\cite{Cortes:1995}.
The problem can be stated as
\begin{equation*}
\begin{array}{ll}
	\mbox{minimize}   & x^\tpose x + \lambda \sum_{i=1}^m \max(0, b_i a_i^T x + 1),
\end{array}
\end{equation*}
where $b_i \in \{ -1, +1 \}$ is a set label, and $a_i$ is a vector of features for the $i$-th point.
The problem can be equivalently represented as the following \gls{QP}
\begin{equation*}
\begin{array}{ll}
    \mbox{minimize}   & x^\tpose x + \lambda \ones^\tpose t \\
    \mbox{subject to}
    & t \ge \diag(b)Ax + \ones \\
    & t \ge 0,
\end{array}
\end{equation*}
where~$\diag(b)$ denotes the diagonal matrix with elements of~$b$ on its diagonal.

\paragraph{Problem instances.}
We choose the vector~$b$ so that
\begin{equation*}
    b_i = \begin{cases}
        +1  & i \le m/2 \\
        -1  & \text{otherwise},
    \end{cases}
\end{equation*}
and the elements of $A$ as
\begin{equation*}
    A_{ij} \sim \begin{cases}
        \mathcal{N}(+1/n,1/n) & i \le m/2 \\
        \mathcal{N}(-1/n,1/n) & \text{otherwise},
    \end{cases}
\end{equation*}
with $15\%$ nonzeros per case.

\begin{figure}[ht!]
\begin{center}
\begin{tikzpicture}
 \begin{groupplot}[
	width=.475\textwidth,
	height=.2\textheight,
	group style={
	group size= 2 by 3,
	vertical sep=3.5em,
	horizontal sep=3.5em},
	grid=both, 
	ymin = 0.0001,
	ymax = 1000,
	xmin = 100,
	xmax = 100000000,
	xmode=log,
	ymode=log
		]
		\nextgroupplot[title=Random QP, ylabel={Computation time}, y unit=\si{\second}]
		\addplot [scatter, only marks, scatter src=explicit symbolic, 
		discard if not={class}{Random QP},
		scatter/classes={Random QP={mark=\gurobimarker, gurobi}}] %
			table[x=N, y=run_time, meta=class, col sep=comma] {data/results/benchmark_problems/GUROBI/results.csv};\label{plot:osqp:timing_gurobi_low}
		\addplot [scatter, only marks, scatter src=explicit symbolic, 
		discard if not={class}{Random QP},
		scatter/classes={Random QP={mark=\osqpmarker, osqp}}] %
			table[x=N, y=run_time, meta=class, col sep=comma] {data/results/benchmark_problems/OSQP/results.csv};\label{plot:osqp:timing_osqp_low}
		\coordinate (top) at (rel axis cs:0,1);
		\nextgroupplot[title=Eq QP]
		\addplot [scatter, only marks, scatter src=explicit symbolic, 
		discard if not={class}{Eq QP},
		scatter/classes={Eq QP={mark=\gurobimarker, gurobi}}] %
			table[x=N, y=run_time, meta=class, col sep=comma] {data/results/benchmark_problems/GUROBI/results.csv};
		\addplot [scatter, only marks, scatter src=explicit symbolic, 
		discard if not={class}{Eq QP},
		scatter/classes={Eq QP={mark=\osqpmarker, osqp}}] %
			table[x=N, y=run_time, meta=class, col sep=comma] {data/results/benchmark_problems/OSQP/results.csv};
		\nextgroupplot[title=Portfolio, ylabel={Computation time}, y unit=\si{\second}]
		\addplot [scatter, only marks, scatter src=explicit symbolic, 
		discard if not={class}{Portfolio},
		scatter/classes={Portfolio={mark=\gurobimarker, gurobi}}] %
			table[x=N, y=run_time, meta=class, col sep=comma] {data/results/benchmark_problems/GUROBI/results.csv};
		\addplot [scatter, only marks, scatter src=explicit symbolic, 
		discard if not={class}{Portfolio},
		scatter/classes={Portfolio={mark=\osqpmarker, osqp}}] %
			table[x=N, y=run_time, meta=class, col sep=comma] {data/results/benchmark_problems/OSQP/results.csv};
		\nextgroupplot[title=Lasso]
		\addplot [scatter, only marks, scatter src=explicit symbolic, 
		discard if not={class}{Lasso},
		scatter/classes={Lasso={mark=\gurobimarker, gurobi}}] %
			table[x=N, y=run_time, meta=class, col sep=comma] {data/results/benchmark_problems/GUROBI/results.csv};
		\addplot [scatter, only marks, scatter src=explicit symbolic, 
		discard if not={class}{Lasso},
		scatter/classes={Lasso={mark=\osqpmarker, osqp}}] %
			table[x=N, y=run_time, meta=class, col sep=comma] {data/results/benchmark_problems/OSQP/results.csv};
		\nextgroupplot[title=SVM, ylabel={Computation time}, y unit=\si{\second}, xlabel={Problem dimension $N$}]
		\addplot [scatter, only marks, scatter src=explicit symbolic, 
		discard if not={class}{SVM},
		scatter/classes={SVM={mark=\gurobimarker, gurobi}}] %
			table[x=N, y=run_time, meta=class, col sep=comma] {data/results/benchmark_problems/GUROBI/results.csv};
		\addplot [scatter, only marks, scatter src=explicit symbolic, 
		discard if not={class}{SVM},
		scatter/classes={SVM={mark=\osqpmarker, osqp}}] %
			table[x=N, y=run_time, meta=class, col sep=comma] {data/results/benchmark_problems/OSQP/results.csv};
		\nextgroupplot[title=Huber, xlabel={Problem dimension $N$}]
		\addplot [scatter, only marks, scatter src=explicit symbolic, 
		discard if not={class}{Huber},
		scatter/classes={Huber={mark=\gurobimarker, gurobi}}] %
			table[x=N, y=run_time, meta=class, col sep=comma] {data/results/benchmark_problems/GUROBI/results.csv};
		\addplot [scatter, only marks, scatter src=explicit symbolic, 
		discard if not={class}{Huber},
		scatter/classes={Huber={mark=\osqpmarker, osqp}}] %
			table[x=N, y=run_time, meta=class, col sep=comma] {data/results/benchmark_problems/OSQP/results.csv};
	    \coordinate (bot) at (rel axis cs:1,0);
\end{groupplot}
	\path (top|-current bounding box.north)--
	      coordinate(legendpos)
	      (bot|-current bounding box.north);
	\matrix[
	    matrix of nodes,
	    anchor=south,
	    draw,
	    inner sep=0.2em,
	] at ([yshift=1ex]legendpos) {
	    \ref{plot:osqp:timing_gurobi_low}& GUROBI &[5pt] \ref{plot:osqp:timing_osqp_low}& OSQP\\};
\end{tikzpicture}
\begin{tikzpicture}
		\begin{axis}[
		width=.475\textwidth,
		height=.2\textheight,
		grid=both, 
		ymin = 0.0001,
		ymax = 1000,
		xmin = 100,
		xmax = 100000000,
		title = {Control},
		xlabel={Problem dimension $N$},
		ylabel={Computation time}, y unit=\si{\second},
		xmode=log,
		ymode=log]
		\addplot [scatter, only marks, scatter src=explicit symbolic, 
		discard if not={class}{Control},
		scatter/classes={Control={mark=\gurobimarker, gurobi}}] %
			table[x=N, y=run_time, meta=class, col sep=comma] {data/results/benchmark_problems/GUROBI/results.csv};
		\addplot [scatter, only marks, scatter src=explicit symbolic, 
		discard if not={class}{Control},
		scatter/classes={Control={mark=\osqpmarker, osqp}}] %
			table[x=N, y=run_time, meta=class, col sep=comma] {data/results/benchmark_problems/OSQP/results.csv};
		\end{axis}
	\end{tikzpicture}
\end{center}
\caption{Computation time vs problem dimension for OSQP and GUROBI for low accuracy mode.}
\label{fig:computation_times_low_accuracy}
\end{figure}

\begin{figure}[ht!]
\begin{center}
\begin{tikzpicture}
 \begin{groupplot}[
	width=.475\textwidth,
	height=.2\textheight,
	group style={
	group size= 2 by 4,
	vertical sep=3.5em,
	horizontal sep=3.5em},
	grid=both, 
	ymin = 0.0001,
	ymax = 1000,
	xmin = 100,
	xmax = 100000000,
	xmode=log,
	ymode=log
		]
		\nextgroupplot[title=Random QP, ylabel={Computation time}, y unit=\si{\second}]
		\addplot [scatter, only marks, scatter src=explicit symbolic, 
		discard if not={class}{Random QP},
		scatter/classes={Random QP={mark=\gurobimarker, gurobi}}] %
			table[x=N, y=run_time, meta=class, col sep=comma] {data/results/benchmark_problems_high_accuracy/GUROBI_high/results.csv};\label{plot:osqp:timing_gurobi_high}
		\addplot [scatter, only marks, scatter src=explicit symbolic, 
		discard if not={class}{Random QP},
		scatter/classes={Random QP={mark=\osqpmarker, osqp}}] %
			table[x=N, y=run_time, meta=class, col sep=comma] {data/results/benchmark_problems_high_accuracy/OSQP_high/results.csv};\label{plot:osqp:timing_osqp_high}
		\coordinate (top) at (rel axis cs:0,1);
		\nextgroupplot[title=Eq QP]
		\addplot [scatter, only marks, scatter src=explicit symbolic, 
		discard if not={class}{Eq QP},
		scatter/classes={Eq QP={mark=\gurobimarker, gurobi}}] %
			table[x=N, y=run_time, meta=class, col sep=comma] {data/results/benchmark_problems_high_accuracy/GUROBI_high/results.csv};
		\addplot [scatter, only marks, scatter src=explicit symbolic, 
		discard if not={class}{Eq QP},
		scatter/classes={Eq QP={mark=\osqpmarker, osqp}}] %
			table[x=N, y=run_time, meta=class, col sep=comma] {data/results/benchmark_problems_high_accuracy/OSQP_high/results.csv};
		\nextgroupplot[title=Portfolio, ylabel={Computation time}, y unit=\si{\second}]
		\addplot [scatter, only marks, scatter src=explicit symbolic, 
		discard if not={class}{Portfolio},
		scatter/classes={Portfolio={mark=\gurobimarker, gurobi}}] %
			table[x=N, y=run_time, meta=class, col sep=comma] {data/results/benchmark_problems_high_accuracy/GUROBI_high/results.csv};
		\addplot [scatter, only marks, scatter src=explicit symbolic, 
		discard if not={class}{Portfolio},
		scatter/classes={Portfolio={mark=\osqpmarker, osqp}}] %
			table[x=N, y=run_time, meta=class, col sep=comma] {data/results/benchmark_problems_high_accuracy/OSQP_high/results.csv};
		\nextgroupplot[title=Lasso]
		\addplot [scatter, only marks, scatter src=explicit symbolic, 
		discard if not={class}{Lasso},
		scatter/classes={Lasso={mark=\gurobimarker, gurobi}}] %
			table[x=N, y=run_time, meta=class, col sep=comma] {data/results/benchmark_problems_high_accuracy/GUROBI_high/results.csv};
		\addplot [scatter, only marks, scatter src=explicit symbolic, 
		discard if not={class}{Lasso},
		scatter/classes={Lasso={mark=\osqpmarker, osqp}}] %
			table[x=N, y=run_time, meta=class, col sep=comma] {data/results/benchmark_problems_high_accuracy/OSQP_high/results.csv};
		\nextgroupplot[title=SVM, ylabel={Computation time}, y unit=\si{\second}, xlabel={Problem dimension $N$}]
		\addplot [scatter, only marks, scatter src=explicit symbolic, 
		discard if not={class}{SVM},
		scatter/classes={SVM={mark=\gurobimarker, gurobi}}] %
			table[x=N, y=run_time, meta=class, col sep=comma] {data/results/benchmark_problems_high_accuracy/GUROBI_high/results.csv};
		\addplot [scatter, only marks, scatter src=explicit symbolic, 
		discard if not={class}{SVM},
		scatter/classes={SVM={mark=\osqpmarker, osqp}}] %
			table[x=N, y=run_time, meta=class, col sep=comma] {data/results/benchmark_problems_high_accuracy/OSQP_high/results.csv};
		\nextgroupplot[title=Huber, xlabel={Problem dimension $N$}]
		\addplot [scatter, only marks, scatter src=explicit symbolic, 
		discard if not={class}{Huber},
		scatter/classes={Huber={mark=\gurobimarker, gurobi}}] %
			table[x=N, y=run_time, meta=class, col sep=comma] {data/results/benchmark_problems_high_accuracy/GUROBI_high/results.csv};
		\addplot [scatter, only marks, scatter src=explicit symbolic, 
		discard if not={class}{Huber},
		scatter/classes={Huber={mark=\osqpmarker, osqp}}] %
			table[x=N, y=run_time, meta=class, col sep=comma] {data/results/benchmark_problems_high_accuracy/OSQP_high/results.csv};
	    \coordinate (bot) at (rel axis cs:1,0);
\end{groupplot}
	\path (top|-current bounding box.north)--
	      coordinate(legendpos)
	      (bot|-current bounding box.north);
	\matrix[
	    matrix of nodes,
	    anchor=south,
	    draw,
	    inner sep=0.2em,
	] at ([yshift=1ex]legendpos) {
	    \ref{plot:osqp:timing_gurobi_high}& GUROBI &[5pt] \ref{plot:osqp:timing_osqp_high}& OSQP\\};
\end{tikzpicture}
\begin{tikzpicture}
		\begin{axis}[
		width=.475\textwidth,
		height=.2\textheight,
		grid=both, 
		ymin = 0.0001,
		ymax = 1000,
		xmin = 100,
		xmax = 100000000,
		title = {Control},
		xlabel={Problem dimension $N$},
		ylabel={Computation time}, y unit=\si{\second},
		xmode=log,
		ymode=log]
		\addplot [scatter, only marks, scatter src=explicit symbolic, 
		discard if not={class}{Control},
		scatter/classes={Control={mark=\gurobimarker, gurobi}}] %
			table[x=N, y=run_time, meta=class, col sep=comma] {data/results/benchmark_problems_high_accuracy/GUROBI_high/results.csv};
		\addplot [scatter, only marks, scatter src=explicit symbolic, 
		discard if not={class}{Control},
		scatter/classes={Control={mark=\osqpmarker, osqp}}] %
			table[x=N, y=run_time, meta=class, col sep=comma] {data/results/benchmark_problems_high_accuracy/OSQP_high/results.csv};
		\end{axis}
	\end{tikzpicture}
\end{center}
\caption{Computation time vs problem dimension for OSQP and GUROBI for high accuracy mode.}
\label{fig:computation_times_high_accuracy}
\end{figure}

\begin{figure}[ht]
	\begin{center}
		\begin{tikzpicture}
		\begin{axis}[
		width=0.9\textwidth,
		height=0.33\textwidth,
		grid=both, 
		grid style={line width=.1pt, draw=gray!30},
		major grid style={line width=.2pt,draw=gray!60},
		title = {Low accuracy},
		xlabel={Performance ratio $\tau$},
				ylabel={Ratio of problems solved},
				ymin = 0,
				ymax = 1,
		xmin = 1,
				xmax = 10000,
		xmode=log,
				legend entries={OSQP, GUROBI, MOSEK, ECOS, qpOASES},
				legend cell align=left,
		legend style={at={(.95,.05)},anchor=south east, fill=white, fill opacity=1, draw opacity=1,text opacity=1},
		log ticks with fixed point,
				every axis plot/.append style={very thick}
				]
		\addplot [osqp, \osqplinestyle] table[x=tau, y=OSQP, col sep=comma] {data/results/benchmark_problems/performance_profiles.csv};
		\addplot [gurobi, \gurobilinestyle] table[x=tau, y=GUROBI, col sep=comma] {data/results/benchmark_problems/performance_profiles.csv};
		\addplot [mosek, \moseklinestyle] table[x=tau, y=MOSEK, col sep=comma] {data/results/benchmark_problems/performance_profiles.csv};
				\addplot [ecos] table[x=tau, y=ECOS, col sep=comma] {data/results/benchmark_problems/performance_profiles.csv};
		\addplot [qpoases, \qpoaseslinestyle] table[x=tau, y=qpOASES, col sep=comma] {data/results/benchmark_problems/performance_profiles.csv};
	    \end{axis}
		\end{tikzpicture}

		\vspace{1em}

		\begin{tikzpicture}
		\begin{axis}[
		width=0.9\textwidth,
		height=0.33\textwidth,
		grid=both, 
		grid style={line width=.1pt, draw=gray!30},
		major grid style={line width=.2pt,draw=gray!60},
		title = {High accuracy},
		xlabel={Performance ratio $\tau$},
				ylabel={Ratio of problems solved},
				ymin = 0,
				ymax = 1,
		xmin = 1,
				xmax = 10000,
		xmode=log,
				legend entries={OSQP, GUROBI, MOSEK, ECOS, qpOASES},
				legend cell align=left,
		legend style={at={(.95,.05)},anchor=south east, fill=white, fill opacity=1, draw opacity=1,text opacity=1},
		log ticks with fixed point,
				every axis plot/.append style={very thick}
				]
		\addplot [osqp, \osqplinestyle] table[x=tau, y=OSQP_high, col sep=comma] {data/results/benchmark_problems_high_accuracy/performance_profiles.csv};
		\addplot [gurobi, \gurobilinestyle] table[x=tau, y=GUROBI_high, col sep=comma] {data/results/benchmark_problems_high_accuracy/performance_profiles.csv};
		\addplot [mosek, \moseklinestyle] table[x=tau, y=MOSEK_high, col sep=comma] {data/results/benchmark_problems_high_accuracy/performance_profiles.csv};
				\addplot [ecos] table[x=tau, y=ECOS_high, col sep=comma] {data/results/benchmark_problems_high_accuracy/performance_profiles.csv};
		\addplot [qpoases, \qpoaseslinestyle] table[x=tau, y=qpOASES, col sep=comma] {data/results/benchmark_problems_high_accuracy/performance_profiles.csv};
	    \end{axis}
		\end{tikzpicture}
    \end{center}
	\caption{Benchmark problems comparison with performance profiles.}
	\label{fig:benchmarks_performance_profiles}
\end{figure}

\begin{figure}[ht]
	\begin{center}
		\begin{tikzpicture}
		\begin{axis}[
		width=0.9\textwidth,
		height=0.33\textwidth,
		grid=both, 
		grid style={line width=.1pt, draw=gray!30},
		major grid style={line width=.2pt,draw=gray!60},
		title = {Low accuracy},
		xlabel={Performance ratio $\tau$},
				ylabel={Ratio of problems solved},
				ymin = 0,
				ymax = 1,
		xmin = 1,
				xmax = 10000,
		xmode=log,
				legend entries={OSQP, GUROBI, MOSEK},
				legend cell align=left,
		legend style={at={(.95,.05)},anchor=south east, fill=white, fill opacity=1, draw opacity=1,text opacity=1},
		log ticks with fixed point,
				every axis plot/.append style={very thick}
				]
		\addplot [osqp, \osqplinestyle] table[x=tau, y=OSQP, col sep=comma] {data/results/suitesparse_problems/performance_profiles.csv};
		\addplot [gurobi, \gurobilinestyle] table[x=tau, y=GUROBI, col sep=comma] {data/results/suitesparse_problems/performance_profiles.csv};
		\addplot [mosek, \moseklinestyle] table[x=tau, y=MOSEK, col sep=comma] {data/results/suitesparse_problems/performance_profiles.csv};
	    \end{axis}
		\end{tikzpicture}

		\vspace{1em}

		\begin{tikzpicture}
		\begin{axis}[
		width=0.9\textwidth,
		height=0.33\textwidth,
		grid=both, 
		grid style={line width=.1pt, draw=gray!30},
		major grid style={line width=.2pt,draw=gray!60},
		title = {High accuracy},
		xlabel={Performance ratio $\tau$},
				ylabel={Ratio of problems solved},
				ymin = 0,
				ymax = 1,
		xmin = 1,
				xmax = 10000,
		xmode=log,
				legend entries={OSQP, GUROBI, MOSEK},
				legend cell align=left,
		legend style={at={(.95,.05)},anchor=south east, fill=white, fill opacity=1, draw opacity=1,text opacity=1},
		log ticks with fixed point,
				every axis plot/.append style={very thick}
				]
		\addplot [osqp, \osqplinestyle] table[x=tau, y=OSQP_high, col sep=comma] {data/results/suitesparse_problems_high_accuracy/performance_profiles.csv};
		\addplot [gurobi, \gurobilinestyle] table[x=tau, y=GUROBI_high, col sep=comma] {data/results/suitesparse_problems_high_accuracy/performance_profiles.csv};
		\addplot [mosek, \moseklinestyle] table[x=tau, y=MOSEK_high, col sep=comma] {data/results/suitesparse_problems_high_accuracy/performance_profiles.csv};
	    \end{axis}
		\end{tikzpicture}
    \end{center}
	\caption{SuiteSparse matrix problems comparison with performance profiles.}
	\label{fig:suitesparse_performance_profiles}
\end{figure}

\begin{figure}[ht]
	\begin{center}
		\begin{tikzpicture}
		\begin{axis}[
		width=0.9\textwidth,
		height=0.33\textwidth,
		grid=both, 
		grid style={line width=.1pt, draw=gray!30},
		major grid style={line width=.2pt,draw=gray!60},
		title = {Low accuracy},
		xlabel={Performance ratio $\tau$},
				ylabel={Ratio of problems solved},
				ymin = 0,
				ymax = 1,
		xmin = 1,
				xmax = 10000,
		xmode=log,
				legend entries={OSQP, GUROBI, MOSEK},
				legend cell align=left,
		legend style={at={(.95,.05)},anchor=south east, fill=white, fill opacity=1, draw opacity=1,text opacity=1},
		log ticks with fixed point,
				every axis plot/.append style={very thick}
				]
		\addplot [osqp, \osqplinestyle] table[x=tau, y=OSQP, col sep=comma] {data/results/maros_meszaros_problems/performance_profiles.csv};
		\addplot [gurobi, \gurobilinestyle] table[x=tau, y=GUROBI, col sep=comma] {data/results/maros_meszaros_problems/performance_profiles.csv};
		\addplot [mosek, \moseklinestyle] table[x=tau, y=MOSEK, col sep=comma] {data/results/maros_meszaros_problems/performance_profiles.csv};
	    \end{axis}
		\end{tikzpicture}

		\vspace{1em}

		\begin{tikzpicture}
		\begin{axis}[
		width=0.9\textwidth,
		height=0.33\textwidth,
		grid=both, 
		grid style={line width=.1pt, draw=gray!30},
		major grid style={line width=.2pt,draw=gray!60},
		title = {High accuracy},
		xlabel={Performance ratio $\tau$},
				ylabel={Ratio of problems solved},
				ymin = 0,
				ymax = 1,
		xmin = 1,
				xmax = 10000,
		xmode=log,
				legend entries={OSQP, GUROBI, MOSEK},
				legend cell align=left,
		legend style={at={(.95,.05)},anchor=south east, fill=white, fill opacity=1, draw opacity=1,text opacity=1},
		log ticks with fixed point,
				every axis plot/.append style={very thick}
				]
		\addplot [osqp, \osqplinestyle] table[x=tau, y=OSQP_high, col sep=comma] {data/results/maros_meszaros_problems_high_accuracy/performance_profiles.csv};
		\addplot [gurobi, \gurobilinestyle] table[x=tau, y=GUROBI_high, col sep=comma] {data/results/maros_meszaros_problems_high_accuracy/performance_profiles.csv};
		\addplot [mosek, \moseklinestyle] table[x=tau, y=MOSEK_high, col sep=comma] {data/results/maros_meszaros_problems_high_accuracy/performance_profiles.csv};
	    \end{axis}
		\end{tikzpicture}
    \end{center}
	\caption{Maros-M\'{e}sz\'{a}ros problems comparison with performance profiles.}
	\label{fig:maros_performance_profiles}
\end{figure}